\documentclass[a4paper,12pt]{amsart}

\usepackage{amssymb}
\usepackage{amsmath}
\usepackage{amsthm}
\usepackage{graphics}

\allowdisplaybreaks

\numberwithin{equation}{section}

\theoremstyle{plain}
\newtheorem{thm}{Theorem}[section]

\newtheorem{lem}[thm]{Lemma}

\theoremstyle{definition}

\newtheorem{rem}[thm]{Remark}

\newcommand{\R}{\mathbb{R}}

\newcommand{\Z}{\mathbb{Z}}

\newcommand{\calF}{\mathcal{F}}

\newcommand{\calS}{\mathcal{S}}

\title[The inclusion between Besov spaces and modulation spaces]
{The dilation property of modulation spaces and  their
inclusion relation with Besov spaces}
\author{Mitsuru Sugimoto \and Naohito Tomita}
\date{}

\address{Mitsuru Sugimoto \\
Department of Mathematics \\
Osaka University \\
Machikaneyama 1-16 \\
Toyonaka, Osaka 560-0043, Japan}
\email{sugimoto@math.wani.osaka-u.ac.jp}

\address{Naohito Tomita \\
Department of Mathematics \\
Osaka University \\
Machikaneyama 1-16 \\
Toyonaka, Osaka 560-0043, Japan}
\email{tomita@gaia.math.wani.osaka-u.ac.jp}



\begin{document}
\maketitle
\begin{abstract}
We consider the dilation property of the modulation spaces
$M^{p,q}$.
Let $D_\lambda:f(t)\mapsto f(\lambda t)$ be the dilation operator, and
we consider the behavior of the operator norm
$\|D_\lambda\|_{M^{p,q}\to M^{p,q}}$ with respect to $\lambda$.
Our result determines the best order for it, and as an application,
we establish the optimality of the inclusion relation
between the modulation spaces and Besov space, which was proved by
Toft \cite{Toft}.
\end{abstract}
\section{Introduction}\label{section1}
The modulation spaces $M^{p,q}$ were first introduced by
Feichtinger in
\cite{Feichtinger} and \cite{Feichtinger2}.
The exact definition will be given in the next section, but
the main idea is to consider the decaying property
of a function with respect to the space variable and the
variable of its Fourier transform simultaneously.
That is exactly the heart of the matter of the time-frequency
analysis which is originated in signal analysis or quantum mechanics.
\par
Based on a similar idea, Sj\"ostrand \cite{Sjostrand}
independently introduced a symbol class which assures the $L^2$-boundedness
of corresponding pseudo-differential operators.
In the last decade, the theory of the modulation spaces has been developed,
and its usefulness for the theory of pseudo-differential operators
is getting realized gradually.
Nowadays Sj\"ostrand's symbol class is recognized as a special case of the
modulation spaces by Gr\"ochenig \cite{Grochenig2}.
Gr\"ochenig and Heil \cite{Grochenig-Heil} also used the
modulation spaces, as a powerful tool, to show trace-class results for
pseudo-differential operators.
Consult Gr\"ochenig \cite{Grochenig} for further and detailed
history of this research fields.
\par
Now we are in a situation to start showing fundamental properties
of the modulation spaces, in order to apply them for many other problems.
Actually in Toft's recent work \cite{Toft}, he investigated the mapping
property of convolutions, and showed Young-type results for the modulation
space.
As an application, he showed an inclusion relation between
the modulation spaces and Besov spaces.
We remark that Besov spaces are used in various problems of partial
differential equations, and his result will help us to understand how they
are translated into the terminology of the modulation spaces.
\par
Among many other important properties to be shown,
we focus on the dilation property of the modulation spaces in this article.
Since $M^{2,2}(\R^n)=L^2(\R^n)$,
we have easily
$\|f_\lambda\|_{M^{2,2}}=\lambda^{-n/2}\|f\|_{M^{2,2}}$
by the change of variables $t\mapsto \lambda^{-1}t$,
where $f_{\lambda}(t)=f(\lambda t)$ and $t\in \R^n$.
But it is not clear how $\|f_\lambda\|_{M^{p,q}}$ behaves like with respect to
$\lambda$ except for the case $(p,q)=(2,2)$.
Our objective is to draw the complete picture of the best order of $\lambda$
for every pair of $(p,q)$ (Theorem \ref{1.1}).
\par
We can expect various kinds of applications of this consideration.
In fact, this kind of dilation property is frequently used in the
\lq\lq scaling argument", which is a popular tool to know
the best possible order of the conditions in problems
of partial differential equations.
Actually, in this article, we also show the best possibility of Toft's
inclusion relation mentioned above, as a side product of the main argument
(Theorem \ref{1.2}).
\par
In order to state our main results, we introduce several indexes.
For $1\le p \le \infty$,
we denote the conjugate exponent of $p$ by $p'$
(that is, $1/p+1/p'=1$).
We define subsets of $(1/p,1/q)\in[0,1]\times[0,1]$ in the following way:
\begin{alignat*}{2}
&I_1\,:\,\max(1/p,1/p')\le 1/q, &\qquad
&I_1^*\,:\,\min(1/p,1/p')\ge 1/q,
\\
&I_2\,:\,\max(1/q,1/2)\le 1/p', &\qquad
&I_2^*\,:\,\min(1/q,1/2)\ge 1/p',
\\
&I_3\,:\,\max(1/q,1/2)\le 1/p, &\qquad
&I_3^*\,:\,\min(1/q,1/2)\ge 1/p.
\end{alignat*}
See the following figure:
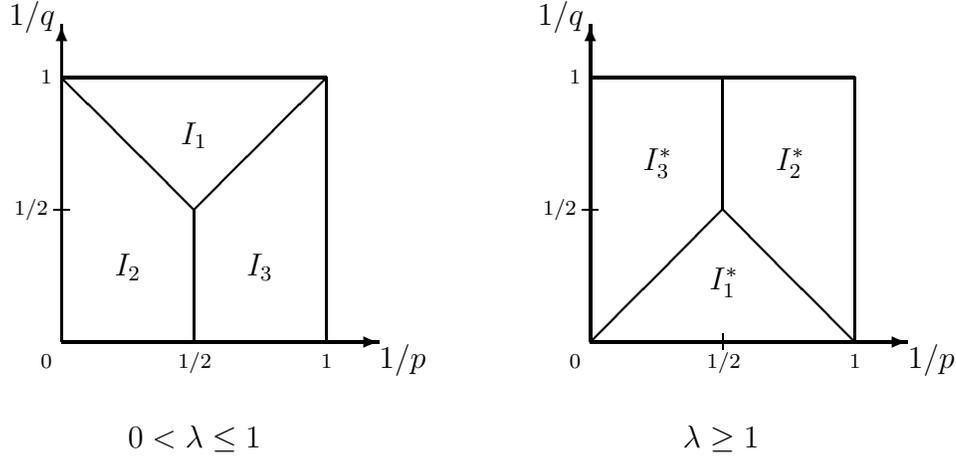
\begin{center}
\begin{picture}(360,190)
\thicklines
\put(45,10){$0<\lambda \le 1$}
\put(20,50){\vector(1,0){120}}
\put(20,150){\line(1,0){100}}
\put(20,50){\vector(0,1){120}}
\put(120,50){\line(0,1){100}}
\put(70,50){\line(0,1){50}}
\put(20,150){\line(1,-1){50}}
\put(70,100){\line(1,1){50}}
\put(65,125){$I_1$}
\put(40,75){$I_2$}
\put(90,75){$I_3$}
\put(140,40){$1/p$}
\put(0,170){$1/q$}
\thinlines
\put(17,100){\line(1,0){6}}
\put(12,40){{\tiny $0$}}
\put(64,40){{\tiny $1/2$}}
\put(118,40){{\tiny $1$}}
\put(2,98){{\tiny $1/2$}}
\put(12,148){{\tiny $1$}}
\thicklines
\put(255,10){$\lambda \ge 1$}
\put(220,50){\vector(1,0){120}}
\put(220,150){\line(1,0){100}}
\put(220,50){\vector(0,1){120}}
\put(320,50){\line(0,1){100}}
\put(270,100){\line(0,1){50}}
\put(220,50){\line(1,1){50}}
\put(270,100){\line(1,-1){50}}
\put(265,70){$I_1^*$}
\put(240,115){$I_3^*$}
\put(290,115){$I_2^*$}
\put(340,40){$1/p$}
\put(200,170){$1/q$}
\thinlines
\put(270,47){\line(0,1){6}}
\put(217,100){\line(1,0){6}}
\put(212,40){{\tiny $0$}}
\put(264,40){{\tiny $1/2$}}
\put(318,40){{\tiny $1$}}
\put(202,98){{\tiny $1/2$}}
\put(212,148){{\tiny $1$}}
\end{picture}
\end{center}
\noindent
In \cite{Toft},
Toft introduced the indexes
\begin{align*}
&\nu_1(p,q)
=\max\{0,1/q-\min(1/p,1/p')\},
\\
&\nu_2(p,q)
=\min\{0,1/q-\max(1/p,1/p')\}.
\end{align*}
Note that
\[
\nu_1(p,q)=
\begin{cases}
0 &\text{if} \quad(1/p,1/q) \in I_1^*, \\
1/p+1/q-1 &\text{if} \quad(1/p,1/q) \in I_2^*, \\
-1/p+1/q &\text{if} \quad (1/p,1/q) \in I_3^*,
\end{cases}
\]
and
\[
\nu_2(p,q)=
\begin{cases}
0 &\text{if} \quad(1/p,1/q) \in I_1, \\
1/p+1/q-1 &\text{if} \quad(1/p,1/q) \in I_2, \\
-1/p+1/q &\text{if} \quad(1/p,1/q) \in I_3.
\end{cases}
\]
We also introduce the indexes
\[
\mu_1(p,q)=\nu_1(p,q)-1/p,
\qquad
\mu_2(p,q)=\nu_2(p,q)-1/p.
\]
Then we have
\[
\mu_1(p,q)=
\begin{cases}
-1/p &\text{if} \quad(1/p,1/q) \in I_1^*, \\
1/q-1 &\text{if} \quad(1/p,1/q) \in I_2^*, \\
-2/p+1/q &\text{if} \quad(1/p,1/q) \in I_3^*,
\end{cases}
\]
and
\[
\mu_2(p,q)=
\begin{cases}
-1/p &\text{if} \quad(1/p,1/q) \in I_1, \\
1/q-1 &\text{if} \quad(1/p,1/q) \in I_2, \\
-2/p+1/q &\text{if} \quad(1/p,1/q) \in I_3.
\end{cases}
\]
\par
Our first main result is on the dilation property of the modulation spaces.
For a function (or tempered distribution)
$f$ on $\R^n$ and $\lambda>0$,
we use the notation $f_{\lambda}$ which is defined by
$f_{\lambda}(t)=f(\lambda t)$, $t\in \R^n$.

\begin{thm}\label{1.1}
Let $1\le p,q \le \infty$.
Then the following are true{\rm :}
\begin{enumerate}
\item
There exists a constant $C>0$ such that
\[
\|f_{\lambda}\|_{M^{p,q}}
\le C\lambda^{n \mu_1(p,q)}\|f\|_{M^{p,q}}
\quad \text{for all} \
f \in M^{p,q}(\R^n)
\ \text{and} \
\lambda \ge 1.
\]
Conversely, if there exist constants
$C>0$ and $\alpha \in \R$ such that
\[
\|f_{\lambda}\|_{M^{p,q}}
\le C\lambda^{\alpha}\|f\|_{M^{p,q}}
\quad \text{for all} \
f \in M^{p,q}(\R^n)
\ \text{and} \
\lambda \ge 1,
\]
then $\alpha \ge n \mu_1(p,q)$.
\item
There exists a constant $C>0$ such that
\[
\|f_{\lambda}\|_{M^{p,q}}
\le C\lambda^{n \mu_2(p,q)}\|f\|_{M^{p,q}}
\quad \text{for all} \
f \in M^{p,q}(\R^n)
\ \text{and} \
0<\lambda \le 1.
\]
Conversely, if there exist constants
$C>0$ and $\beta \in \R$ such that
\[
\|f_{\lambda}\|_{M^{p,q}}
\le C\lambda^{\beta}\|f\|_{M^{p,q}}
\quad \text{for all} \
f \in M^{p,q}(\R^n)
\ \text{and} \
0<\lambda \le 1,
\]
then $\beta \le n \mu_2(p,q)$.
\end{enumerate}
\end{thm}

Since the Gauss function $\varphi(t)=e^{-|t|^2}$ does not change its form
under the Fourier transformation, 
the modulation norm of it can have a \lq\lq good" property.
In this sense, it is reasonable to believe that the Gauss function
$f=\varphi$ attains the critical order of $\|f_\lambda\|_{M^{p,q}}$
with respect to $\lambda$.
But it is not true because
$\|\varphi_{\lambda}\|_{M^{p,q}}\sim \lambda^{n(1/q-1)}$
in the case $\lambda \ge 1$
and $\|\varphi_{\lambda}\|_{M^{p,q}}\sim \lambda^{-n/p}$
in the case $0<\lambda \le 1$ (see Lemma \ref{2.1}).
Theorem \ref{1.1} says that they are not critical orders for every pair of
$(p,q)$.
\par
It should be pointed out here that
the behavior of $\|f_\lambda\|_{M^{p,q}}$ with respect to $\lambda$
might depend on the choice of $f\in M^{p,q}(\R^n)$.
In fact,
$f(t)=\sum_{k \in \Z^n}e^{ik\cdot t}\, \psi(t-k)$,
where $\psi$ is an appropriate Schwartz function,
has the property
$\|f_{\lambda}\|_{M^{p,\infty}}\sim \lambda^{-2n/p}$ ($0<\lambda \le 1$)
in the case $1 \le p \le 2$ (Lemma \ref{3.9}),
while the Gauss function has the different behavior
$\|\varphi_{\lambda}\|_{M^{p,\infty}}\sim \lambda^{-n/p}$
($0<\lambda \le 1$) as mentioned above.
On the other hand, the $L^p$-norm never has such a property
since $\|f_\lambda\|_{L^p}=\lambda^{-n/p}\|f\|_{L^p}$
for all $f\in L^p(\R^n)$.
That is one of great differences between the modulation spaces
and $L^p$-spaces.
\par
Our second main result is on the optimality of the inclusion
relation between the modulation spaces and Besov spaces.
In \cite[Theorem 3.1]{Toft}, Toft proved the inclusions
\[
B_{n \nu_1(p,q)}^{p,q}(\R^n)
\hookrightarrow M^{p,q}(\R^n)
\hookrightarrow B_{n \nu_2(p,q)}^{p,q}(\R^n)
\]
for $1\le p,q \le \infty$.
He also remarked that the left inclusion
is optimal in the case $1\le p=q\le2$,
that is, if $B_{s_1}^{p,p}(\R^n) \hookrightarrow M^{p,p}(\R^n)$
then $s_1 \ge n\nu_1(p,p)$.
The same is true for the right inclusion in the case $2\le p=q\le\infty$,
that is, if $M^{p,p}(\R^n) \hookrightarrow B_{s_2}^{p,p}(\R^n)$
then $s_2 \le n\nu_2(p,p)$
(\cite[Remark 3.11]{Toft}).
The next theorem says that Toft's inclusion result
is optimal in the above meaning for every pair of $(p,q)$.

\begin{thm}\label{1.2}
Let $1\le p,q \le \infty$ and $s \in \R$.
Then the following are true{\rm :}
\begin{enumerate}
\item
If $B_s^{p,q}(\R^n) \hookrightarrow M^{p,q}(\R^n)$,
then $s \ge n \nu_1(p,q)$.
\item
If $M^{p,q}(\R^n) \hookrightarrow B_s^{p,q}(\R^n)$
and $1\le p,q < \infty$,
then $s \le n \nu_2(p,q)$.
\end{enumerate}
\end{thm}

We end this introduction by explaining the plan of this article.
In Section 2, we give the precise definition and basic properties
of the modulation spaces and Besov spaces.
In Sections 3 and 4, we prove Theorems \ref{1.1} and \ref{1.2}
respectively.
\section{Preliminaries}\label{section2}
We introduce the modulation spaces
based on Gr\"ochenig \cite{Grochenig}.
Let $\calS(\R^n)$ and $\calS'(\R^n)$ be the Schwartz spaces of
rapidly decreasing smooth functions
and tempered distributions,
respectively.
We define the Fourier transform $\hat{f}$
and the inverse Fourier transform $\calF^{-1}f$
of $f \in \calS(\R^n)$ by
\[
\hat{f}(\xi)
=\int_{\R^n}e^{-i\xi \cdot x}\, f(x)\, dx
\]
and
\[
\calF^{-1}f(x)
=\frac{1}{(2\pi)^n}
\int_{\R^n}e^{ix\cdot \xi}\, f(\xi)\, d\xi.
\]
Fix a function $\varphi \in \calS(\R^n)\setminus \{ 0 \}$
(called the window function).
Then the short-time Fourier transform $V_{\varphi}f$ of
$f \in \calS'(\R^n)$ with respect to $\varphi$
is defined by
\[
V_{\varphi}f(x,\xi)
=\langle f, M_{\xi}T_x\varphi \rangle
\quad \text{for} \ x, \xi \in \R^n,
\]
where $M_{\xi}T_x\varphi(t)=e^{i\xi \cdot t}\varphi(t-x)$
and $\langle \cdot, \cdot \rangle$ is the inner product on $L^2(\R^n)$.
We can express it in a form of the integral
\[
V_{\varphi}f(x,\xi)
=\int_{\R^n}f(t)\, \overline{\varphi(t-x)}\,
e^{-i\xi \cdot t}\, dt,
\]
which has actually the meaning
for an appropriate function $f$ on $\R^n$.
We note that,
for $f \in \calS'(\R^n)$,
$V_{\varphi}f$ is continuous on $\R^{2n}$
and $|V_{\varphi}f(x,\xi)|\le C(1+|x|+|\xi|)^N$
for some constants $C,N \ge 0$
(\cite[Theorem 11.2.3]{Grochenig}).
Let $1\le p,q \le \infty$.
Then the modulation space $M^{p,q}(\R^n)$
consists of all $f \in \calS'(\R^n)$
such that
\[
\|f\|_{M^{p,q}}
=\|V_{\varphi}f\|_{L^{p,q}}
=\left\{ \int_{\R^n} \left(
\int_{\R^n} |V_{\varphi}f(x,\xi)|^{p}\, dx
\right)^{q/p} d\xi \right\}^{1/q}
< \infty.
\]
We note that
$M^{2,2}(\R^n)=L^2(\R^n)$
(\cite[Proposition 11.3.1]{Grochenig})
and $M^{p,q}(\R^n)$ is a Banach space
(\cite[Proposition 11.3.5]{Grochenig}).
The definition of $M^{p,q}(\R^n)$ is independent
of the choice of the window function
$\varphi \in \calS(\R^n)\setminus \{ 0 \}$,
that is,
different window functions
yield equivalent norms
(\cite[Proposition 11.3.2]{Grochenig}).

We also introduce Besov spaces.
Let $1\le p,q \le \infty$ and $s \in \R$.
Suppose that $\eta ,\psi \in \calS(\R^n)$ satisfy
${\rm supp}\, \eta \subset \{\xi: |\xi|\le 2\}$,
${\rm supp}\, \psi \subset \{\xi: 1/2 \le |\xi| \le 2\}$
and $\eta(\xi)+\sum_{j=1}^{\infty}\psi(\xi/2^j)=1$
for all $\xi \in \R^n$.
Set $\varphi_0=\eta$
and $\varphi_j=\psi(\cdot/2^j)$
if $j \ge 1$.
Then Besov space $B_s^{p,q}(\R^n)$
consists of all $f \in \calS'(\R^n)$
such that
\[
\|f\|_{B_s^{p,q}}
=\left( \sum_{j=0}^{\infty}
2^{jsq}\|\calF^{-1}[\varphi_j \, \hat{f}]\|_{L^p}^q
\right)^{1/q}
=\left( \sum_{j=0}^{\infty}
2^{jsq}\|\Phi_j * f\|_{L^p}^q
\right)^{1/q}<\infty,
\]
where $\Phi_j=\calF^{-1}\varphi_j$.
We remark
$B^{p,q}_s(\R^n)^*=B^{p',q'}_{-s}(\R^n)$
for $1\le p,q<\infty$.

Finally,
we list below the lemmas
which will be used in the subsequent section.
In this article,
we frequently use the function
$\varphi(t)=e^{-|t|^2}$
which is called the Gauss function.
\begin{lem}[{\cite[Lemma 1.8]{Toft}}]\label{2.1}
Let $\varphi$ be the Gauss function.
Then
\[
\|V_{\varphi}(\varphi_{\lambda})\|_{L^{p,q}}
=\pi^{n(1/p+1/q+1)/2}p^{-n/2p}q^{-n/2q}2^{n/q}
\lambda^{-n/p}(1+\lambda^2)^{n(1/p+1/q-1)/2}.
\]
\end{lem}
Lemma \ref{2.1} says that
$\|\varphi_{\lambda}\|_{M^{p,q}}\sim \lambda^{n(1/q-1)}$
in the case $\lambda \ge 1$
and $\|\varphi_{\lambda}\|_{M^{p,q}}\sim \lambda^{-n/p}$
in the case $0<\lambda \le 1$.
\begin{lem}[{\cite[Corollary 11.2.7]{Grochenig}}]\label{2.2}
Let $f \in \calS'(\R^n)$ and $\varphi,\psi,\gamma \in \calS(\R^n)$.
Then
\[
\langle f,\varphi \rangle
=\frac{1}{\langle \gamma, \psi \rangle}
\int_{\R^{2n}}
V_{\psi}f(x,\xi)\, \overline{V_{\gamma}\varphi(x,\xi)}\,
dx \, d\xi
\quad
\text{for all} \ \varphi \in \calS(\R^n).
\]
\end{lem}
\begin{lem}[{\cite[Lemma 11.3.3]{Grochenig}}]\label{2.3}
Let $f \in \calS'(\R^n)$ and
$\varphi, \psi, \gamma \in \calS(\R^n)$.
Then
\[
|V_{\varphi}f(x,\xi)|
\le \frac{1}{|\langle \gamma, \psi \rangle |}
(|V_{\psi}f|*|V_{\varphi}\gamma|)(x,\xi)
\quad \text{for all} \ x, \xi \in \R^n.
\]
\end{lem}
\begin{lem}[{\cite[Proposition 11.3.4 and Theorem 11.3.6]{Grochenig}}]
\label{2.4}
If $1 \le p,q <\infty$,
then $\calS(\R^n)$ is dense in $M^{p,q}(\R^n)$
and $M^{p,q}(\R^n)^*=M^{p',q'}(\R^n)$
under the duality
\[
\langle f,g \rangle_M
=\frac{1}{\|\varphi\|_{L^2}^2}
\int_{\R^{2n}}
V_{\varphi}f(x,\xi)\, \overline{V_{\varphi}g(x,\xi)}\,
dx \, d\xi
\]
for $f \in M^{p,q}(\R^n)$ and $g \in M^{p',q'}(\R^n)$.
\end{lem}
By Lemmas \ref{2.2} and \ref{2.4},
if $1<p,q \le \infty$ and $f \in M^{p,q}(\R^n)$
then
\begin{equation}\label{(2.1)}
\|f\|_{M^{p,q}}
=\sup \left| \langle f,g \rangle_M \right|
=\sup \left| \langle f,g \rangle \right|
\end{equation}
where the supremum is taken over all
$g \in \calS(\R^n)$ such that
$\|g\|_{M^{p',q'}}=1$.
\begin{lem}[{\cite[Corollary 2.3]{Feichtinger}}]\label{2.5}
Let $1\le p_1,p_2,q_1,q_2 \le \infty$ and $q_2<\infty$.
If $T$ is a linear operator such that
\[
\|Tf\|_{M^{p_1,q_1}}\le
A_1 \|f\|_{M^{p_1,q_1}}
\quad \text{for all} \ f \in M^{p_1,q_1}(\R^n)
\]
and
\[
\|Tf\|_{M^{p_2,q_2}}\le
A_2 \|f\|_{M^{p_2,q_2}}
\quad \text{for all} \ f \in M^{p_2,q_2}(\R^n),
\]
then
\[
\|Tf\|_{M^{p,q}}\le
CA_1^{1-\theta}A_2^{\theta}\|f\|_{M^{p,q}}
\quad \text{for all} \ f \in M^{p,q}(\R^n),
\]
where $1/p=(1-\theta)/p_1+\theta/p_2$,
$1/q=(1-\theta)/q_1+\theta/q_2$,
$0 \le \theta \le 1$
and $C$ is independent of $T$.
\end{lem}
\begin{rem}\label{2.6}
Lemma \ref{2.5} with the case $q_1=q_2=\infty$
is treated in \cite[Remark 3.2]{Toft},
which says that it is true under a modification.
\end{rem}
\section{The dilation property of modulation spaces}\label{section3}
In this section,
we prove Theorem \ref{1.1}
which appeared in the introduction.
We begin by preparing
the following lemma:
\begin{lem}\label{3.1}
Let $1\le p,q \le \infty$.
Then there exists a constant $C>0$
such that
\[
\|f_{\lambda}\|_{M^{p,q}}
\le C\lambda^{-n(1/p-1/q+1)}(1+\lambda^2)^{n/2}
\|f\|_{M^{p,q}}
\]
for all $f \in M^{p,q}(\R^n)$
and $\lambda >0$.
\end{lem}
\begin{proof}
Let $\varphi$ be the Gauss function,
that is, $\varphi(t)=e^{-|t|^2}$.
By a change of variable,
we have
\[
\|f_{\lambda}\|_{M^{p,q}}
=\|V_{\varphi}(f_{\lambda})\|_{L^{p,q}}
=\lambda^{-n(1/p-1/q+1)}\|V_{\varphi_{1/\lambda}}f\|_{L^{p,q}}.
\]
From Lemma \ref{2.3} it follows that
\[
\left|V_{\varphi_{1/\lambda}}f(x,\xi) \right|
\le \|\varphi\|_{L^2}^{-2}
(|V_{\varphi}f|*|V_{\varphi_{1/\lambda}}\varphi|)(x,\xi).
\]
Hence, by Young's inequality and Lemma \ref{2.1},
we get
\begin{align*}
\|f_{\lambda}\|_{M^{p,q}}
&\le \lambda^{-n(1/p-1/q+1)}
\|\varphi\|_{L^2}^{-2}
\|V_{\varphi_{1/\lambda}}\varphi\|_{L^{1,1}}
\|V_{\varphi}f\|_{L^{p,q}} \\
&=\lambda^{-n(1/p-1/q+1)}
\|\varphi\|_{L^2}^{-2}
\|V_{\varphi}(\varphi_{1/\lambda})\|_{L^{1,1}}
\|V_{\varphi}f\|_{L^{p,q}} \\
&=\lambda^{-n(1/p-1/q+1)}\|\varphi\|_{L^2}^{-2}
\left( \pi^{3n/2}2^n
(\lambda^{-1})^{-n}(1+\lambda^{-2})^{n/2} \right)
\|f\|_{M^{p,q}} \\
&=C_n\lambda^{-n(1/p-1/q+1)}
(1+\lambda^{2})^{n/2} \|f\|_{M^{p,q}}.
\end{align*}
The proof is complete.
\end{proof}
We are now ready to prove
Theorem \ref{1.1} (1) with $(1/p,1/q) \in I_1^*$
and (2) with $(1/p,1/q) \in I_1$.

\medskip
\noindent
{\it Proof of Theorem \ref{1.1} (2) with $(1/p,1/q) \in I_1$}.
Suppose that $(1/p,1/q) \in I_1$.
Then $\mu_2(p,q)=-1/p$.
Let $1\le r \le \infty$.
By Lemma \ref{3.1},
we have
\begin{equation}\label{(3.1)}
\|f_{\lambda}\|_{M^{r,1}}
\le C\lambda^{-n/r}\|f\|_{M^{r,1}}
\quad \text{for all} \
f \in M^{r,1}(\R^n)
\ \text{and} \
0<\lambda \le 1.
\end{equation}
On the other hand,
since $M^{2,2}(\R^n)=L^2(\R^n)$,
we have
\begin{equation}\label{(3.2)}
\|f_{\lambda}\|_{M^{2,2}}
\le C\lambda^{-n/2}\|f\|_{M^{2,2}}
\quad \text{for all} \
f \in M^{2,2}(\R^n)
\ \text{and} \
0<\lambda \le 1.
\end{equation}
Take $1 \le r \le \infty$
and $0 \le \theta \le 1$
such that $1/p=(1-\theta)/r+\theta/2$
and $1/q=(1-\theta)/1+\theta/2$.
Then, by the interpolation theorem
(Lemma \ref{2.5}),
(\ref{(3.1)}) and (\ref{(3.2)}) give
\[
\|f_{\lambda}\|_{M^{p,q}}
\le C\left( \lambda^{-n/r}\right)^{1-\theta}
\left( \lambda^{-n/2} \right)^{\theta}
\|f\|_{M^{p,q}}
\]
for all $f \in M^{p,q}(\R^n)$
and
$0<\lambda \le 1$.
Since $(1-\theta)/r=1/p+1/q-1$ and $\theta/2=1-1/q$,
we get
\begin{equation}\label{(3.3)}
\|f_{\lambda}\|_{M^{p,q}}
\le C\lambda^{-n/p}
\|f\|_{M^{p,q}}
\quad \text{for all} \
f \in M^{p,q}(\R^n)
\ \text{and} \
0<\lambda \le 1.
\end{equation}
This is the first part of Theorem \ref{1.1} (2)
with $(1/p,1/q) \in I_1$.

We next prove the second part of Theorem \ref{1.1} (2)
with $(1/p,1/q) \in I_1$.
Let $(1/p,1/q) \in I_1$.
Assume that there exist constants
$C>0$ and $\beta \in \R$ such that
\[
\|f_{\lambda}\|_{M^{p,q}}
\le C\lambda^{\beta}
\|f\|_{M^{p,q}}
\quad \text{for all} \
f \in M^{p,q}(\R^n)
\ \text{and} \
0<\lambda \le 1.
\]
Let $\varphi$ be the Gauss function.
We note that the Gauss function belongs to $M^{p,q}(\R^n)$.
Then, by Lemma \ref{2.1} and our assumption,
we have
\begin{align*}
C_{p,q}\lambda^{-n/p}
&\le C_{p,q}\lambda^{-n/p}
(1+\lambda^2)^{n(1/p+1/q-1)/2} \\
&=\|V_{\varphi}(\varphi_{\lambda})\|_{L^{p,q}}
=\|\varphi_{\lambda}\|_{M^{p,q}}
\le C\lambda^{\beta}\|\varphi\|_{M^{p,q}}
\end{align*}
for all $0<\lambda \le 1$.
This is possible only if $\beta \le -n/p$.
The proof is complete.

\medskip
\noindent
{\it Proof of Theorem \ref{1.1} (1) with $(1/p,1/q) \in I_1^*$}.
We note that $\mu_1(p,q)=-1/p$ if $(1/p,1/q) \in I_1^*$.
Let $1\le p \le \infty$ and $q \ge 2$
be such that $(1/p,1/q) \in I_1^*$.
Then $(1/p',1/q') \in I_1$.
We first consider the case $p \neq 1$.
Since $1< p,q \le \infty$,
by duality
(\ref{(2.1)})
and Theorem \ref{1.1} (2) with $(1/p',1/q') \in I_1$,
we have
\begin{align*}
\|f_{\lambda}\|_{M^{p,q}}
&=\sup \left| \langle f_{\lambda},g \rangle \right|
=\lambda^{-n}\sup \left| \langle f,g_{1/\lambda} \rangle \right| \\
&\le \lambda^{-n}\sup \|f\|_{M^{p,q}}\|g_{1/\lambda}\|_{M^{p',q'}} \\
&\le \lambda^{-n}\sup \|f\|_{M^{p,q}}
\left( C(\lambda^{-1})^{-n/p'}\|g\|_{M^{p',q'}} \right)
=C\lambda^{-n/p}\|f\|_{M^{p,q}}
\end{align*}
for all $f \in M^{p,q}(\R^n)$
and $\lambda \ge 1$,
where the supremum is taken over all $g \in \calS(\R^n)$
such that $\|g\|_{M^{p',q'}}=1$.
In the case $p=1$,
by Lemma \ref{3.1},
we see that
\[
\|f_{\lambda}\|_{M^{1,\infty}}
\le C\lambda^{-n}\|f\|_{M^{1,\infty}}
\quad \text{for all} \
f \in M^{1,\infty}(\R^n)
\ \text{and} \
\lambda \ge 1.
\]
Hence,
we obtain the first part of Theorem \ref{1.1} (1)
with $(1/p,1/q) \in I_1^*$.

We consider the second part of Theorem \ref{1.1} (1)
with $(1/p,1/q) \in I_1^*$.
Let $1 < p < \infty$ and $2 \le q < \infty$
be such that
$(1/p,1/q) \in I_1^*$.
Assume that there exist constants
$C>0$ and $\alpha \in \R$ such that
\[
\|g_{\lambda}\|_{M^{p,q}}
\le C\lambda^{\alpha}
\|g\|_{M^{p,q}}
\quad \text{for all} \
g \in M^{p,q}(\R^n)
\ \text{and} \
\lambda \ge 1.
\]
Then, by duality and our assumption,
we have
\begin{align*}
\|f_{\lambda}\|_{M^{p',q'}}
&=\sup \left| \langle f_{\lambda},g \rangle \right|
=\lambda^{-n}\sup \left| \langle f,g_{1/\lambda} \rangle \right| \\
&\le \lambda^{-n}
\sup \|f\|_{M^{p',q'}}\|g_{1/\lambda}\|_{M^{p,q}} \\
&\le \lambda^{-n}
\sup \|f\|_{M^{p',q'}}
\left( C(\lambda^{-1})^{\alpha}\|g\|_{M^{p,q}}\right)
=C\lambda^{-n-\alpha}\|f\|_{M^{p',q'}}
\end{align*}
for all $f \in M^{p',q'}(\R^n)$
and $0<\lambda \le 1$,
where the supremum is taken over all $g \in \calS(\R^n)$
such that $\|g\|_{M^{p,q}}=1$.
Since $(1/p',1/q') \in I_1$,
by Theorem \ref{1.1} (2) with $(1/p',1/q') \in I_1$,
we get $-n-\alpha \le -n/p'$.
This implies $\alpha \ge -n/p$.

We next consider the case $q=\infty$.
Let $1 \le r \le \infty$.
Assume that there exist constants
$C>0$ and $\alpha \in \R$ such that
\begin{equation}\label{(3.4)}
\|f_{\lambda}\|_{M^{r,\infty}}
\le C\lambda^{\alpha}
\|f\|_{M^{r,\infty}}
\quad \text{for all} \
f \in M^{r,\infty}(\R^n)
\ \text{and} \
\lambda \ge 1,
\end{equation}
where $\alpha <-n/r$.
Since $M^{2,2}(\R^n)=L^2(\R^n)$,
we have
\begin{equation}\label{(3.5)}
\|f_{\lambda}\|_{M^{2,2}}
\le C\lambda^{-n/2}\|f\|_{M^{2,2}}
\quad \text{for all} \
f \in M^{2,2}(\R^n)
\ \text{and} \
\lambda \ge 1.
\end{equation}
Then, by the interpolation theorem,
(\ref{(3.4)}) and (\ref{(3.5)}) give
\[
\|f_{\lambda}\|_{M^{p,q}}
\le C\left( \lambda^{\alpha}\right)^{1-\theta}
\left( \lambda^{-n/2}\right)^{\theta}
\|f\|_{M^{p,q}}
\]
for all $f \in M^{p,q}(\R^n)$
and $\lambda \ge 1$,
where $1/p=(1-\theta)/r+\theta/2$,
$1/q=(1-\theta)/\infty+\theta/2$,
and $0 < \theta <1$.
Since $0<\theta <1$,
we note that $1<p<\infty$,
$2<q<\infty$ and $(1/p,1/q) \in I_1^*$.
Since $p=q$ if $r=\infty$,
using that
$1-\theta=r(1/p-1/q)$ if $1\le r <\infty$,
$1-\theta=1-2/q$ if $r=\infty$
and ${\theta}/2=1/q$,
we have
\begin{align*}
\|f_{\lambda}\|_{M^{p,q}}&\le
\begin{cases}
C\lambda^{\alpha r(1/p-1/q)-n/q}
\|f\|_{M^{p,q}}, &\text{if} \quad 1\le r <\infty \\
C\lambda^{\alpha(1-2/q)-n/q}
\|f\|_{M^{p,q}}, &\text{if} \quad  r=\infty
\end{cases} \\
&=\begin{cases}
C\lambda^{\alpha r(1/p-1/q)-n/q+n/p-n/p}
\|f\|_{M^{p,q}}, &\text{if} \quad 1\le r <\infty \\
C\lambda^{\alpha(1-2/q)-n/p}
\|f\|_{M^{p,q}}, &\text{if} \quad r=\infty
\end{cases} \\
&=\begin{cases}
C\lambda^{(\alpha r+n)(1/p-1/q)-n/p}
\|f\|_{M^{p,q}}, &\text{if} \quad 1\le r <\infty \\
C\lambda^{\alpha(1-2/q)-n/p}
\|f\|_{M^{p,q}}, &\text{if} \quad r=\infty
\end{cases}
\end{align*}
for all $f \in M^{p,q}(\R^n)$
and $\lambda \ge 1$.
However,
since $(\alpha r+n)(1/p-1/q)<0$ if $1\le r <\infty$
and $\alpha(1-2/q)<0$ if $r=\infty$,
this contradicts
Theorem \ref{1.1} (1) with $1<p<\infty$,
$2<q<\infty$ and $(1/p,1/q) \in I_1^*$.
Therefore, $\alpha$ must satisfy $\alpha \ge -n/r$.
The proof is complete.

\medskip
Our next goal is to prove
Theorem \ref{1.1} (1) with $(1/p,1/q) \in I_2^*$
and (2) with $(1/p,1/q) \in I_2$.
\begin{lem}\label{3.2}
Let $1 \le p,q \le \infty$ be such that
$(1/p,1/q) \in I_2^*$
and $1/p \ge 1/q$.
Then there exists a constant $C>0$
such that
\[
\|f_{\lambda}\|_{M^{p,q}}
\le C\lambda^{-n(2/p-1/q)}
(1+\lambda^2)^{n(1/p-1/2)}
\|f\|_{M^{p,q}}
\]
for all $f \in M^{p,q}(\R^n)$
and $\lambda >0$.
\end{lem}
\begin{proof}
Let $1 \le r \le \infty$.
By Lemma \ref{3.1},
we have
\begin{equation}\label{(3.6)}
\|f_{\lambda}\|_{M^{1,r}}
\le C\lambda^{n(1/r-2)}
(1+\lambda^2)^{n/2}\|f\|_{M^{1,r}}
\end{equation}
for all
$f \in M^{1,r}(\R^n)$
and
$\lambda > 0$.
Let $1 \le p,q \le \infty$ be such that
$(1/p,1/q) \in I_2^*$ and $1/p \ge 1/q$.
Take $1\le r \le \infty$ and $0 \le \theta \le 1$
such that
$1/p=(1-\theta)/1+\theta/2$
and $1/q=(1-\theta)/r+\theta/2$.
Then, by the interpolation theorem,
(\ref{(3.2)}), (\ref{(3.5)})
and (\ref{(3.6)})
give
\[
\|f_{\lambda}\|_{M^{p,q}}
\le C\left( \lambda^{n(1/r-2)}
(1+\lambda^2)^{n/2} \right)^{1-\theta}
\left( \lambda^{-n/2}\right)^{\theta}
\|f\|_{M^{p,q}}
\]
for all $f \in M^{p,q}(\R^n)$
and $\lambda >0$.
Using $(1-\theta)/r=1/p+1/q-1$,
$1-\theta=2/p-1$
and ${\theta}/2=-1/p+1$,
we get
\begin{align*}
\|f_{\lambda}\|_{M^{p,q}}
&\le C\lambda^{n\left( (1-\theta)/r-2(1-\theta)-\theta/2 \right)}
(1+\lambda^2)^{n(1-\theta)/2}
\|f\|_{M^{p,q}} \\
&= C\lambda^{-n(2/p-1/q)}
(1+\lambda^2)^{n(1/p-1/2)}
\|f\|_{M^{p,q}}
\end{align*}
for all $f \in M^{p,q}(\R^n)$
and $\lambda >0$.
The proof is complete.
\end{proof}
The proof of the following lemma is
based on that of \cite[Theorem 3]{Triebel}.
\begin{lem}\label{3.3}
Suppose that
$\varphi \in \calS(\R^n)$ is a real-valued function
satisfying
$\varphi=1$ on $[-1/2,1/2]^n$,
${\rm supp}\, \varphi \subset [-1,1]^n$,
$\varphi(t)=\varphi(-t)$
and
$\sum_{k \in \Z^n}\varphi (t-k)=1$
for all $t \in \R^n$.
Then
\[
\sup_{k \in \Z^n}\|(M_k \Phi)*f\|_{L^2}
\le \|V_{\Phi}f\|_{L^{2,\infty}}
\le 5^n\|\Phi\|_{L^1}\sup_{k \in \Z^n}\|(M_k \Phi)*f\|_{L^2}
\]
for all $f \in M^{2,\infty}(\R^n)$,
where $\Phi=\calF^{-1}\varphi$ and $M_k\Phi(t)=e^{ik\cdot t}\Phi(t)$.
\end{lem}
\begin{proof}
Let $f \in M^{2,\infty}(\R^n)$.
Since $\Phi$ is a real-valued function
and $\Phi(t)=\Phi(-t)$ for all $t$,
we have
\begin{align}\label{(3.7)}
|V_{\Phi}f(x,\xi)|
&=\left| \int_{\R^n}f(t)\, \overline{\Phi(t-x)}\,
e^{-i\xi \cdot t}\, dt \right| \\ 
&=\left| \int_{\R^n}f(t)\, \Phi(x-t)\,
e^{i\xi \cdot (x-t)}\, dt \right|
=\left| (M_{\xi}\Phi)*f(x) \right|. \nonumber
\end{align}
We first prove
\begin{equation}\label{(3.8)}
{\rm ess \, sup}_{\xi \in \R^n}
\left(\int_{\R^n}|V_{\Phi}f(x,\xi)|^2 dx\right)^{1/2}
=\sup_{\xi \in \R^n}
\left(\int_{\R^n}|V_{\Phi}f(x,\xi)|^2 dx\right)^{1/2}.
\end{equation}
To prove (\ref{(3.8)}),
it is enough to show that
$\left(\int_{\R^n}|V_{\Phi}f(x,\xi)|^2 dx\right)^{1/2}$
is continuous with respect to $\xi$.
Since
${\rm ess \, sup}_{\xi \in \R^n}
\left(\int_{\R^n}|V_{\Phi}f(x,\xi)|^2 dx\right)^{1/2}<\infty,
$
for each $k \in \Z^n$
there exists $\xi_k \in k/2+[-1/4,1/4]^n$
such that
$\left(\int_{\R^n}|V_{\Phi}f(x,\xi_k)|^2 dx\right)^{1/2}<\infty$.
Then, by (\ref{(3.7)}),
we have
\[
\frac{1}{(2\pi)^{n/2}}
\|\varphi(\cdot -\xi_k)\, \hat{f}\|_{L^2}
=\|(M_{\xi_k}\Phi)*f\|_{L^2}
=\left(\int_{\R^n}|V_{\Phi}f(x,\xi_k)|^2 dx\right)^{1/2}<\infty.
\]
Since
$k/2+[-1/4,1/4]^n \subset \xi_k+[-1/2,1/2]^n$
and  $\varphi(\cdot-\xi_k)=1$
on $\xi_k+[-1/2,1/2]^n$,
we see that $|\hat{f}|^2$ is integrable on $k/2+[-1/4,1/4]^n$.
The arbitrariness of $k \in \Z^n$
gives $\hat{f} \in L_{{\rm loc}}^2(\R^n)$.
By the Lebesgue dominated convergence theorem,
we see that
$\|\varphi(\cdot-\xi)\, \hat{f}\|_{L^2}$
is continuous with respect to $\xi$.
Hence,
$\left(\int_{\R^n}|V_{\Phi}f(x,\xi)|^2 dx\right)^{1/2}$
is continuous with respect to $\xi$.
We obtain (\ref{(3.8)}).
Then, from (\ref{(3.7)}) and (\ref{(3.8)})
it follows that
\begin{align*}
\sup_{k \in \Z^n}\|(M_k\Phi)*f\|_{L^2}
&\le \sup_{\xi \in \R^n}\|(M_{\xi}\Phi)*f\|_{L^2}
=\sup_{\xi \in \R^n}
\left(\int_{\R^n}|V_{\Phi}f(x,\xi)|^2 dx\right)^{1/2} \\
&={\rm ess \, sup}_{\xi \in \R^n}
\left(\int_{\R^n}|V_{\Phi}f(x,\xi)|^2 dx\right)^{1/2}
=\|V_{\Phi}f\|_{L^{2,\infty}}.
\end{align*}

We next prove
$\|V_{\Phi}f\|_{L^{2,\infty}} \le
(5^n\|\Phi\|_{L^1}) \sup_{k \in \Z^n}\|(M_k \Phi)*f\|_{L^2}$.
Let $\xi \in \R^n$.
Since
\begin{align*}
M_{\xi}\Phi
&=\calF^{-1}[\varphi(\cdot-\xi)]
=\calF^{-1}\left[ \varphi(\cdot-\xi)\,
\left( \sum_{k \in \Z^n}\varphi(\cdot-k)\right) \right] \\
&=\sum_{\scriptstyle |k_j-\xi_j|\le 2, \atop \scriptstyle j=1,\cdots , n}
\calF^{-1}[\varphi(\cdot-\xi)\, \varphi(\cdot-k)]
=\sum_{\scriptstyle |k_j-\xi_j|\le 2, \atop \scriptstyle j=1,\cdots , n}
(M_{\xi}\Phi)*(M_{k}\Phi),
\end{align*}
by (\ref{(3.7)}),
we have
\[
|V_{\Phi}f(x,\xi)|
=\left| (M_{\xi}\Phi)*f(x) \right|
\le \sum_{\scriptstyle |k_j-\xi_j|\le 2, \atop \scriptstyle j=1,\cdots , n}
\left| (M_{\xi}\Phi)*(M_{k}\Phi)*f(x) \right|.
\]
Hence, by (\ref{(3.8)}),
we get
\begin{align*}
\|V_{\Phi}f\|_{L^{2,\infty}}
&=\sup_{\xi \in \R^n}
\left( \int_{\R^n}|V_{\Phi}f(x,\xi)|^2 \, dx \right)^{1/2} \\
&\le \sup_{\xi \in \R^n}
\sum_{\scriptstyle |k_j-\xi_j|\le 2, \atop \scriptstyle j=1,\cdots , n}
\| (M_{\xi}\Phi)*(M_{k}\Phi)*f \|_{L^2} \\
&\le \sup_{\xi \in \R^n}
\sum_{\scriptstyle |k_j-\xi_j|\le 2, \atop \scriptstyle j=1,\cdots , n}
\| M_{\xi}\Phi \|_{L^1}\, \|(M_{k}\Phi)*f \|_{L^2} \\
&\le \|\Phi \|_{L^1}
\left( \sup_{\ell \in \Z^n}\|(M_{\ell}\Phi)*f \|_{L^2} \right)
\left( \sup_{\xi \in \R^n}
\sum_{\scriptstyle |k_j-\xi_j|\le 2, \atop \scriptstyle j=1,\cdots , n} 1
\right) \\
&\le 5^n \|\Phi \|_{L^1}
\sup_{\ell \in \Z^n}\|(M_{\ell}\Phi)*f \|_{L^2}.
\end{align*}
The proof is complete.
\end{proof}
We remark that
Lemma \ref{3.1} implies
\[
\|f_{\lambda}\|_{M^{2,\infty}}
\le C\lambda^{-3n/2} \|f\|_{M^{2,\infty}}
\quad \text{for all} \ f \in M^{2,\infty}(\R^n)
\ \text{and} \ 0<\lambda \le 1.
\]
This is not our desired order of $\lambda$
in the case $(p,q)=(2,\infty)$.
But we have
\begin{lem}\label{3.4}
There exists a constant $C>0$ such that
\[
\|f_{\lambda}\|_{M^{2,\infty}}
\le C\lambda^{-n} \|f\|_{M^{2,\infty}}
\quad \text{for all} \ f \in M^{2,\infty}(\R^n)
\ \text{and} \ 0<\lambda \le 1.
\]
\end{lem}
\begin{proof}
Let $\Phi=\calF^{-1}\varphi$,
where $\varphi$ is as in Lemma \ref{3.3}.
Suppose that $f \in M^{2,\infty}(\R^n)$.
We note that $\hat{f} \in L_{\rm loc}^2(\R^n)$
(see the proof of Lemma \ref{3.3}).
Then, by Lemma \ref{3.3},
we see that
\begin{align*}
\|V_{\Phi}(f_{\lambda})\|_{M^{2,\infty}}
&\le 5^n \|\Phi\|_{L^1}
\sup_{k \in \Z^n}\|(M_k\Phi)*f_{\lambda}\|_{L^2} \\
&=(2\pi)^{-n/2} 5^n \|\Phi\|_{L^1}
\sup_{k \in \Z^n}\|\varphi(\cdot-k)\,
\widehat{f_{\lambda}}\|_{L^2} \\
&=C_n \lambda^{-n/2}\sup_{k \in \Z^n}
\|\varphi(\lambda \cdot-k)\, \hat{f}\|_{L^2} \\
&=C_n \lambda^{-n/2}\sup_{k \in \Z^n}
\left\| \varphi(\lambda \cdot-k)
\left( \sum_{\ell \in \Z^n}\varphi(\cdot-\ell)\right)
\hat{f}\right\|_{L^2}.
\end{align*}
Since
\begin{align*}
\left| \varphi(\lambda t-k)
\left( \sum_{\ell \in \Z^n}\varphi(t-\ell)\right)
\hat{f}(t)\right|^2
&\le 4^n \sum_{\ell \in \Z^n}
\left| \varphi(\lambda t-k)\,
\varphi(t-\ell)\,
\hat{f}(t)\right|^2 \\
&=4^n \sum_{\scriptstyle |\ell_j-k_j/\lambda | \le 2/\lambda,
\atop \scriptstyle j=1,\cdots, n}
\left| \varphi(\lambda t-k)\,
\varphi(t-\ell)\,
\hat{f}(t)\right|^2,
\end{align*}
we have
\begin{align*}
&\left\| \varphi(\lambda \cdot-k)
\left( \sum_{\ell \in \Z^n}\varphi(\cdot-\ell)\right)
\hat{f}\right\|_{L^2} \\
&\le \left(
4^n \sum_{\scriptstyle |\ell_j-k_j/\lambda | \le 2/\lambda,
\atop \scriptstyle j=1,\cdots, n}
\int_{\R^n}
\left| \varphi(\lambda t-k)\,
\varphi(t-\ell)\,
\hat{f}(t)\right|^2
\, dt \right)^{1/2} \\
&\le \left(
4^n (2\pi)^n\|\varphi \|_{L^{\infty}}^2
\sum_{\scriptstyle |\ell_j-k_j/\lambda | \le 2/\lambda,
\atop \scriptstyle j=1,\cdots, n}
\|(M_{\ell}\Phi)*f\|_{L^2}^2
\right)^{1/2} \\
&\le \left(
4^n (2\pi)^n\|\varphi \|_{L^{\infty}}^2
\left( \sup_{m \in \Z^n}\|(M_{m}\Phi)*f\|_{L^2} \right)^2
\sum_{\scriptstyle |\ell_j-k_j/\lambda | \le 2/\lambda,
\atop \scriptstyle j=1,\cdots, n} 1
\right)^{1/2} \\
&\le \left(
C_n \|\varphi \|_{L^{\infty}}^2 \lambda^{-n} 
\left( \sup_{m \in \Z^n}\|(M_{m}\Phi)*f\|_{L^2} \right)^2
\right)^{1/2} \\
&=C_n \|\varphi\|_{L^{\infty}}
\lambda^{-n/2}\sup_{m \in \Z^n}\|(M_{m}\Phi)*f\|_{L^2}.
\end{align*}
Hence,
by Lemma \ref{3.3},
we get
\[
\|f_{\lambda}\|_{M^{2,\infty}}
\le C_n
\lambda^{-n}\sup_{m \in \Z^n}\|(M_{m}\Phi)*f\|_{L^2}
\le C_n
\lambda^{-n} \|f\|_{M^{2,\infty}}.
\]
The proof is complete.
\end{proof}
\begin{lem}\label{3.5}
Let $1 \le p \le \infty$. Then the following are true{\rm :}
\begin{enumerate}
\item
If $p \le 2$,
then there exists a constant $C>0$ such that
\[
\|f_{\lambda}\|_{M^{p,1}}
\le C\|f\|_{M^{p,1}}
\quad \text{for all} \ f \in M^{p,1}(\R^n)
\ \text{and} \ \lambda \ge 1.
\]
\item
If $p \ge 2$,
then there exists a constant $C>0$ such that
\[
\|f_{\lambda}\|_{M^{p,1}}
\le C\lambda^{-n(2/p-1)} \|f\|_{M^{p,1}}
\quad \text{for all} \ f \in M^{p,1}(\R^n)
\ \text{and} \ \lambda \ge 1.
\]
\end{enumerate}
\end{lem}
\begin{proof}
We first consider the case $p \le 2$.
By Lemmas \ref{2.2} and \ref{2.4},
and Lemma \ref{3.4},
we have
\begin{align*}
\|f_{\lambda}\|_{M^{2,1}}
&=\sup \left| \langle f_{\lambda},g \rangle_M \right|
=\sup \left| \langle f_{\lambda},g \rangle \right| \\
&=\lambda^{-n}
\sup \left| \langle f,g_{1/\lambda} \rangle \right|
\le \lambda^{-n}
\sup \|f\|_{M^{2,1}}\|g_{1/\lambda}\|_{M^{2,\infty}} \\
&\le \lambda^{-n}
\sup \|f\|_{M^{2,1}}
\left( C\|g\|_{M^{2,\infty}}(1/\lambda)^{-n}\right)
=C\|f\|_{M^{2,1}}
\end{align*}
for all $f \in \calS(\R^n)$ and $\lambda \ge 1$,
where the supremum is taken over all
$g \in M^{2,\infty}(\R^n)$
such that $\|g\|_{M^{2,\infty}}(\R^n)=1$.
Since $\calS(\R^n)$ is dense in $M^{2,1}(\R^n)$,
this gives
\begin{equation}\label{(3.9)}
\|f_{\lambda}\|_{M^{2,1}}
\le C\|f\|_{M^{2,1}}
\quad \text{for all} \ f \in M^{2,1}(\R^n)
\ \text{and} \ \lambda \ge 1.
\end{equation}
On the other hand,
by Lemma \ref{3.1},
we see that
\begin{equation}\label{(3.10)}
\|f_{\lambda}\|_{M^{1,1}}
\le C\|f\|_{M^{1,1}}
\quad \text{for all} \ f \in M^{1,1}(\R^n)
\ \text{and} \ \lambda \ge 1.
\end{equation}
Hence, by the interpolation theorem,
(\ref{(3.9)}) and (\ref{(3.10)}) give
Lemma \ref{3.5} (1).

We next consider the case $p \ge 2$.
By Lemma \ref{3.1},
we have
\begin{equation}\label{(3.11)}
\|f_{\lambda}\|_{M^{\infty,1}}
\le C\lambda^n \|f\|_{M^{\infty,1}}
\quad \text{for all} \
f \in M^{\infty,1}(\R^n)
\ \text{and} \
\lambda \ge 1.
\end{equation}
Therefore, by the interpolation theorem,
(\ref{(3.9)}) and (\ref{(3.11)}) give
\[
\|f_{\lambda}\|_{M^{p,1}}
\le C\left( \lambda^{0}\right)^{1-\theta}
\left( \lambda^{n} \right)^{\theta}
\|f\|_{M^{p,1}}
\]
for all $f \in M^{p,1}(\R^n)$
and
$\lambda \ge 1$,
where 
$1/p=(1-\theta)/2+\theta/\infty$
and $0 \le \theta \le 1$.
Since $\theta=-2/p+1$,
this implies Lemma \ref{3.5} (2).
The proof is complete.
\end{proof}
\begin{lem}\label{3.6}
Let $1 \le p \le \infty$. Then the following are true{\rm :}
\begin{enumerate}
\item
If $p \le 2$,
then there exists a constant $C>0$ such that
\[
\|f_{\lambda}\|_{M^{p,\infty}}
\le C\lambda^{-2n/p}\|f\|_{M^{p,\infty}}
\quad \text{for all} \ f \in M^{p,\infty}(\R^n)
\ \text{and} \ 0<\lambda \le 1.
\]
\item
If $p \ge 2$,
then there exists a constant $C>0$ such that
\[
\|f_{\lambda}\|_{M^{p,\infty}}
\le C\lambda^{-n} \|f\|_{M^{p,\infty}}
\quad \text{for all} \ f \in M^{p,\infty}(\R^n)
\ \text{and} \ 0<\lambda \le 1.
\]
\end{enumerate}
\end{lem}
\begin{proof}
Let $1<p\le 2$.
By duality and Lemma \ref{3.5} (2),
we have
\begin{align*}
\|f_{\lambda}\|_{M^{p,\infty}}
&=\sup \left| \langle f_{\lambda},g \rangle \right|
=\lambda^{-n}
\sup \left| \langle f,g_{1/\lambda} \rangle \right| \\
&\le \lambda^{-n}\sup \|f\|_{M^{p,\infty}}
\left( C(1/\lambda)^{-n(2/p'-1)}\|g\|_{M^{p',1}} \right)
=C\lambda^{-2n/p}\|f\|_{M^{p,\infty}}
\end{align*}
for all $\ f \in M^{p,\infty}(\R^n)$
and $0<\lambda \le 1$,
where the supremum is taken over all
$g \in \calS(\R^n)$ such that $\|g\|_{M^{p',1}}=1$.
In the case $p=1$,
by Lemma \ref{3.1},
we have
\[
\|f_{\lambda}\|_{M^{1,\infty}}
\le C\lambda^{-2n} \|f\|_{M^{1,\infty}}
\quad \text{for all} \ f \in M^{1,\infty}(\R^n)
\ \text{and} \ 0<\lambda \le 1.
\]
Hence, we obtain Lemma \ref{3.6} (1).
In the same way,
we can prove Lemma \ref{3.6} (2).
\end{proof}
\begin{lem}\label{3.7}
Let $1\le p,q \le \infty$,
$(p,q)\neq (1,\infty), (\infty,1)$
and $\epsilon>0$.
Set
\[
f(t)=\sum_{k \neq 0}
|k|^{-n/q-\epsilon}\, e^{ik\cdot t}\, \varphi(t)
\quad \text{in} \quad \calS'(\R^n),
\]
where $\varphi$ is the Gauss function.
Then $f \in M^{p,q}(\R^n)$
and there exists a constant $C>0$
such that
$\|f_{\lambda}\|_{M^{p,q}}\ge C\lambda^{n(1/q-1)+\epsilon}$
for all $0<\lambda \le 1$.
\end{lem}
\begin{proof}
We first prove $f \in M^{p,q}(\R^n)$.
Since
\begin{align*}
&\left| \int_{\R^n}e^{ik\cdot t}\, \varphi(t)\,
\varphi(t-x)\, e^{-i\xi \cdot t}\, dt \right| \\
&=\left| \int_{\R^n}\, \varphi(t)\,
\varphi(x-t) \left\{(1+|\xi-k|^2)^{-n}(I-\Delta_t)^n\,
e^{-i(\xi-k) \cdot t}\right\} dt \right| \\
&=(1+|\xi-k|^2)^{-n}
\left| \sum_{|\beta_1+\beta_2| \le 2n}
C_{\beta_1,\beta_2}
\int_{\R^n}\, (\partial^{\beta_1}\varphi)(t)\,
(\partial^{\beta_2}\varphi)(x-t)\,
e^{-i(\xi-k) \cdot t}\, dt \right| \\
&\le C(1+|\xi-k|^2)^{-n}
\sum_{|\beta_1+\beta_2| \le 2n}
|\partial^{\beta_1}\varphi|*
|\partial^{\beta_2}\varphi|(x),
\end{align*}
we have
\begin{align*}
&\|f\|_{M^{p,q}}
=\| V_{\varphi}f \|_{L^{p,q}} \\
&=\left\{ \int_{\R^n}\left(
\int_{\R^n}\left| \sum_{k \neq 0}|k|^{-n/q-\epsilon}
\int_{\R^n}e^{ik\cdot t}\, \varphi(t)\,
\varphi(t-x)\, e^{-i\xi \cdot t}\, dt \right| dx
\right)^{q/p} d\xi \right\}^{1/q} \\
&\le C\left\{ \int_{\R^n} \left(
\sum_{k \neq 0}|k|^{-n/q-\epsilon}\, (1+|\xi-k|^2)^{-n}
\right)^q d\xi \right\}^{1/q} \\
&=C\left\{ \sum_{\ell \in \Z^n}
\int_{\ell +[-1/2,1/2]^n} \left(
\sum_{k \neq 0}|k|^{-n/q-\epsilon}\, (1+|\xi-k|^2)^{-n}
\right)^q d\xi \right\}^{1/q} \\
&\le C\left\{ \sum_{\ell \in \Z^n}
\left( \sum_{k \neq 0}
|k|^{-n/q-\epsilon}\, (1+|\ell-k|^2)^{-n}
\right)^q \right\}^{1/q}.
\end{align*}
Since $\{|k|^{-n/q-\epsilon}\}_{k \neq 0} \in \ell^q(\Z^n)$,
by Young's inequality,
we see that $f \in M^{p,q}(\R^n)$.

We next consider the second part.
Since $\varphi \in M^{p',q'}(\R^n)$,
by duality,
we have
\begin{align*}
\|f_{\lambda}\|_{M^{p,q}}
&=\sup_{\|g\|_{M^{p',q'}}=1}
\left| \langle f_{\lambda},g \rangle \right|
\ge \left| \langle f_{\lambda},\varphi \rangle \right| \\
&=\left| \sum_{k \neq 0}
|k|^{-n/q-\epsilon}
\int_{\R^n}e^{i(\lambda k)\cdot t}\,
\varphi(\lambda t)\,
\varphi(t)\, dt \right| \\
&=\pi^{n/2}(1+\lambda^2)^{-n/2}
\sum_{k \neq 0}
|k|^{-n/q-\epsilon}\,
e^{-\frac{\lambda^2 |k|^2}{4(1+\lambda^2)}} \\
&\ge (\pi/2)^{n/2}
\sum_{\scriptstyle 0<|k_j|\le 1/\lambda,
\atop \scriptstyle j=1,\cdots,n}
|k|^{-n/q-\epsilon}\,
e^{-\frac{\lambda^2 |k|^2}{4(1+\lambda^2)}} \\
&\ge C\lambda^{n/q+\epsilon}
\sum_{\scriptstyle 0<|k_j|\le 1/\lambda,
\atop \scriptstyle j=1,\cdots,n}1
\ge C\lambda^{n(1/q-1)+\epsilon}
\end{align*}
for all $0<\lambda \le 1$.
The proof is complete.
\end{proof}
We are now ready to prove Theorem \ref{1.1} (1)
with $(1/p,1/q) \in I_2^*$
and (2) with $(1/p,1/q) \in I_2$.

\medskip
\noindent
{\it Proof of Theorem \ref{1.1} (2) with $(1/p,1/q) \in I_2$}.
We note that $\mu_2(p,q)=1/q-1$ if $(1/p,1/q) \in I_2$.
Let $p \ge 2$ and $1 \le q \le \infty$ be such that
$(1/p,1/q) \in I_2$ and $1/p \le 1/q$.
If $q=1$ then $p=\infty$,
and we have already proved this case
in Theorem \ref{1.1} (2)
with $(1/p,1/q) \in I_1$.
Hence, we may assume $1<q\le \infty$.
We note that $1 \le p' \le 2$ and $1\le q'<\infty$.
Since $(1/p', 1/q') \in I_2^*$
and $1/p' \ge 1/q'$,
by duality and Lemma \ref{3.2},
we have
\begin{align}\label{(3.12)}
&\|f_{\lambda}\|_{M^{p,q}}
=\sup \left| \langle f_{\lambda},g \rangle \right|
=\lambda^{-n}
\sup \left| \langle f,g_{1/\lambda} \rangle \right| \\
&\le \lambda^{-n}
\sup \|f\|_{M^{p,q}}
\|g_{1/\lambda}\|_{M^{p',q'}} \nonumber \\
&\le \lambda^{-n}
\sup \|f\|_{M^{p,q}}
\left( C(\lambda^{-1})^{-n(2/p'-1/q')}
(1+\lambda^{-2})^{n(1/p'-1/2)}
\|g\|_{M^{p',q'}} \right) \nonumber \\
&\le C\lambda^{n(1/q-1)}\|f\|_{M^{p,q}} \nonumber
\end{align}
for all $f \in M^{p,q}(\R^n)$
and $0<\lambda \le 1$,
where the supremum is taken over all
$g \in \calS(\R^n)$
such that $\|g\|_{M^{p',q'}}=1$.
This is the first part of Theorem \ref{1.1} (2)
with $(1/p,1/q) \in I_2$
and $1/p \le 1/q$.
Let $p \ge 2$ and $2\le q <\infty$ be such that
$(1/p,1/q) \in I_2$
and $1/p \ge 1/q$.
From (\ref{(3.12)}) it follows that
\begin{equation}\label{(3.13)}
\|f_{\lambda}\|_{M^{r,r}}
\le C\lambda^{n(1/r-1)}\|f\|_{M^{r,r}}
\quad \text{for all} \
f \in M^{r,r}(\R^n)
\ \text{and} \
0<\lambda \le 1,
\end{equation}
where $2 \le r \le \infty$.
Take $2 \le r \le \infty$ and $0 \le \theta \le 1$
such that
$1/p=(1-\theta)/2+{\theta}/r$
and $1/q=(1-\theta)/\infty+{\theta}/r$.
Since $q \neq \infty$,
we note that $r \neq \infty$.
Then, by the interpolation theorem,
Lemma \ref{3.4} and (\ref{(3.13)})
give
\[
\|f_{\lambda}\|_{M^{p,q}}
\le C\left( \lambda^{-n} \right)^{1-\theta}
\left( \lambda^{n(1/r-1)} \right)^{\theta}
\|f\|_{M^{p,q}}
\]
for all $f \in M^{p,q}(\R^n)$
and $0<\lambda \le 1$.
Since ${\theta}/r=1/q$,
we have
\[
\|f_{\lambda}\|_{M^{p,q}}
\le C\lambda^{n({\theta}/r-1)}\|f\|_{M^{p,q}}
\le C\lambda^{n(1/q-1)}\|f\|_{M^{p,q}}
\]
for all $f \in M^{p,q}(\R^n)$
and $0<\lambda \le 1$.
In the case $q=\infty$,
by Lemma \ref{3.6} (2),
we have nothing to prove.
Hence, we obtain the first part of Theorem \ref{1.1} (2)
with $(1/p,1/q) \in I_2$
and $1/p \ge 1/q$.

We next consider the second part of Theorem \ref{1.1} (2)
with $(1/p,1/q) \in I_2$.
Let $p \ge 2$ and $1 \le q \le \infty$ be such that
$(1/p,1/q) \in I_2$.
Since $(1/\infty,1/1) \in I_1$,
we may assume $(p,q)\neq (\infty,1)$.
Assume that there exist constants
$C>0$ and $\beta \in \R$ such that
\[
\|f_{\lambda}\|_{M^{p,q}}
\le C\lambda^{\beta}\|f\|_{M^{p,q}}
\quad \text{for all} \
f \in M^{p,q}(\R^n)
\ \text{and} \
0<\lambda \le 1,
\]
where $\beta > n(1/q-1)$.
Then we can take $\epsilon>0$
such that $n(1/q-1)+\epsilon<\beta$.
For this $\epsilon$,
we set
\[
f(t)=\sum_{k \neq 0}
|k|^{-n/q-\epsilon}\, e^{ik\cdot t}\, \varphi(t),
\]
where $\varphi$ is the Gauss function.
Then,
by Lemma \ref{3.7},
we see that
$f \in M^{p,q}(\R^n)$ and
there exists a constant $C'>0$
such that
$\|f_{\lambda}\|_{M^{p,q}}\ge C'\lambda^{n(1/q-1)+\epsilon}$
for all $0<\lambda \le 1$.
Hence,
\[
C'\lambda^{n(1/q-1)+\epsilon}
\le \|f_{\lambda}\|_{M^{p,q}}
\le C\lambda^{\beta}\|f\|_{M^{p,q}}
\]
for all $0<\lambda \le 1$.
However,
since $n(1/q-1)+\epsilon<\beta$,
this is contradiction.
Therefore, $\beta$ must satisfy $\beta \le n(1/q-1)$.
The proof is complete.

\medskip
\noindent
{\it Proof of Theorem \ref{1.1} (1) with $(1/p,1/q) \in I_2^*$}.
We note that $\mu_1(p,q)=1/q-1$ if $(1/p,1/q) \in I_2^*$.
In every case except for $(p,q)\neq (1,\infty)$,
by duality,
Theorem \ref{1.1} (2)
with $(1/p',1/q') \in I_2$
and the same argument as in the proof of Theorem \ref{1.1} (1)
with $(1/p,1/q) \in I_1^*$,
we can prove Theorem \ref{1.1} (1)
with $(1/p,1/q) \in I_2^*$.
For the case $(p,q)=(1,\infty)$,
we have already proved in Theorem \ref{1.1} (1)
with $(1/p,1/q) \in I_1^*$.

\medskip
Our last goal of this section is to prove
Theorem \ref{1.1} (1) with $(1/p,1/q) \in I_3^*$
and (2) with $(1/p,1/q) \in I_3$.
In the following lemma,
we use the fact
that there exists $\varphi \in \calS(\R^n)$
such that ${\rm supp}\, \varphi \subset [-1/8,1/8]^n$
and $|\hat{\varphi}|\ge 1$ on $[-2,2]^n$
(see, for example,
the proof of \cite[Theorem 2.6]{Frazier-Jawerth}).
\begin{lem}\label{3.8}
Let $1 \le p \le \infty$,
$1 \le q <\infty$ and $\epsilon>0$.
Suppose that $\varphi, \psi \in \calS(\R^n)$ satisfy
${\rm supp}\, \varphi \subset [-1/8,1/8]^n$,
${\rm supp}\, \psi \subset [-1/2,1/2]^n$,
$|\hat{\varphi}|\ge 1$ on $[-2,2]^n$
and $\psi=1$ on $[-1/4,1/4]^n$.
Set
\[f(t)=\sum_{k \neq 0}
|k|^{-n/q-\epsilon}\,
e^{ik\cdot t}\,
\psi(t-k)
\quad \text{in} \quad \calS'(\R^n).
\]
Then $f \in M^{p,q}(\R^n)$
and there exists a constant $C>0$
such that
\[
\|V_{\varphi}(f_{\lambda})\|_{L^{p,q}}
\ge C\lambda^{-n(2/p-1/q)+\epsilon}
\quad
\text{for all} \ 0<\lambda \le 1.
\]
\end{lem}
\begin{proof}
In the same way as in the proof of Lemma \ref{3.7},
we can prove $f \in M^{p,q}(\R^n)$.
We consider the second part.
Since
$\|V_{\varphi}(f_{\lambda})\|_{L^{p,q}}
=\lambda^{-n(1/p-1/q+1)}\|V_{\varphi_{1/\lambda}}f\|_{L^{p,q}}$,
it is enough to show that
$\|V_{\varphi_{1/\lambda}}f\|_{L^{p,q}}
\ge C\lambda^{-n/p+n+\epsilon}$
for all $0<\lambda \le 1$.
We note that
${\rm supp}\, \varphi((\cdot-x)/\lambda) \subset \ell + [-1/4,1/4]^n$
for all $0<\lambda \le 1, \ell \in \Z^n$
and $x \in \ell+[-1/8,1,8]^n$.
Since
${\rm supp}\, \psi(\cdot-k) \subset k+[-1/2,1/2]^n$
and $\psi(t-k)=1$ if $t \in k+[-1/4, 1/4]^n$,
it follows that
\begin{align*}
&\left( \int_{\R^n}\left| V_{\varphi_{1/\lambda}}f(x,\xi)
\right|^p\, dx \right)^{1/p} \\
&=\left( \int_{\R^n}\left|
\sum_{k \neq 0}|k|^{-n/q-\epsilon}
\int_{\R^n}e^{ik\cdot t}\, \psi(t-k)\,
\overline{\varphi\left(\frac{t-x}{\lambda}\right)}\,
e^{-i\xi \cdot t}\, dt
\right|^p\, dx \right)^{1/p} \\
&\ge \left( \sum_{\ell \neq 0} \int_{\ell+[-1/8,1/8]^n}
\left| \sum_{k \neq 0}|k|^{-n/q-\epsilon}
\int_{\R^n}e^{-i(\xi-k) \cdot t}\, \psi(t-k)\,
\overline{\varphi\left(\frac{t-x}{\lambda}\right)}\, dt
\right|^p\, dx \right)^{1/p} \\
&=\left( \sum_{\ell \neq 0} \int_{\ell+[-1/8,1/8]^n}
\left| |\ell|^{-n/q-\epsilon}
\int_{\R^n}e^{-i(\xi-\ell) \cdot t}\,
\overline{\varphi\left(\frac{t-x}{\lambda}\right)}\, dt
\right|^p\, dx \right)^{1/p} \\
&=4^{-n/p}\left( \sum_{\ell \neq 0}
\left| |\ell|^{-n/q-\epsilon} \lambda^n
\hat{\varphi}(-\lambda(\xi-\ell))
\right|^p \right)^{1/p}.
\end{align*}
Hence,
using $|\hat{\varphi}|\ge 1$ on $[-2,2]^n$,
we get
\begin{align*}
\|V_{\varphi_{1/\lambda}}f\|_{L^{p,q}}
&\ge 4^{-n/p}\lambda^n
\left\{ \int_{\R^n}
\left( \sum_{\ell \neq 0}
\left| |\ell|^{-n/q-\epsilon}
\hat{\varphi}(-\lambda(\xi-\ell))
\right|^p \right)^{q/p} d\xi \right\}^{1/q} \\
&=4^{-n/p}\lambda^{n-n/q}
\left\{ \int_{\R^n}
\left( \sum_{\ell \neq 0}
\left| |\ell|^{-n/q-\epsilon}
\hat{\varphi}(\xi+\lambda \ell))
\right|^p \right)^{q/p} d\xi \right\}^{1/q} \\
&\ge 4^{-n/p}\lambda^{n-n/q}
\left\{ \int_{[-1,1]^n}
\left( \sum_{\scriptstyle 0<|\ell_j|\le 1/\lambda,
\atop \scriptstyle j=1,\cdots,n}
\left| |\ell|^{-n/q-\epsilon}
\hat{\varphi}(\xi+\lambda \ell))
\right|^p \right)^{q/p} d\xi \right\}^{1/q} \\
&\ge 4^{-n/p}2^{n/q}\lambda^{n-n/q}
\left( \sum_{\scriptstyle 0<|\ell_j|\le 1/\lambda,
\atop \scriptstyle j=1,\cdots,n}
|\ell|^{-(n/q+\epsilon)p} \right)^{1/p} \\
&\ge C_n \lambda^{n-n/q}
\lambda^{n/q+\epsilon}
\left( \sum_{\scriptstyle 0<|\ell_j|\le 1/\lambda,
\atop \scriptstyle j=1,\cdots,n}1 \right)^{1/p}
\ge C_n \lambda^{-n/p+n+\epsilon}
\end{align*}
for all $0<\lambda \le 1$.
The proof is complete.
\end{proof}
For Lemma \ref{3.8},
we do not need $\epsilon>0$ in the case $q=\infty$.
\begin{lem}\label{3.9}
Let $1 \le p \le \infty$.
Suppose that $\varphi,\psi \in \calS(\R^n)$
are as in Lemma \ref{3.8}.
Set
\[
f(t)=
\sum_{k \in \Z^n}e^{ik\cdot t}\, \psi(t-k)
\quad \text{in} \quad \calS'(\R^n).
\]
Then $f \in M^{p,\infty}(\R^n)$
and there exists a constant $C>0$
such that $\|f_{\lambda}\|_{M^{p,\infty}} \ge C \lambda^{-2n/p}$
for all $0<\lambda \le 1$.
In particular, if $1 \le p \le 2$
then there exist constants
$C,C'>0$ such that
\[
C\lambda^{-2n/p}\le 
\|f_{\lambda}\|_{M^{p,\infty}}\le C'\lambda^{-2n/p}
\quad \text{for all} \ 0<\lambda \le 1.
\]
\end{lem}
\begin{proof}
In the same way as in the proof of Lemma \ref{3.7},
we can prove
\[
\left| \int_{\R^n}e^{ik\cdot t}\,
\psi(t-k)\, \overline{\varphi(t-x)}\,
e^{-i\xi \cdot t}\, dt \right|
\le C (1+|x-k|^2)^{-n}\, (1+|\xi-k|^2)^{-n}.
\]
Hence,
\begin{align*}
\left| V_{\varphi}f(x,\xi)\right|
&=\left| \sum_{k \in \Z^n}
\int_{\R^n}e^{ik\cdot t}\,
\psi(t-k)\, \overline{\varphi(t-x)}\,
e^{-i\xi \cdot t}\, dt \right| \\
&\le C\sum_{k \in \Z^n}
(1+|x-k|^2)^{-n}\, (1+|\xi-k|^2)^{-n}
\le C(1+|x-\xi|^2)^{-n}
\end{align*}
for all $x, \xi \in \R^n$.
This implies $f \in M^{p,\infty}(\R^n)$.

We next consider the second part.
Since
$\|V_{\varphi_{1/\lambda}}f(\cdot,\xi)\|_{L^p}$
is continuous with respect to $\xi \in \R^n$,
we see that
$\|V_{\varphi_{1/\lambda}}f\|_{L^{p,\infty}}
=\sup_{\xi \in \R^n}\|V_{\varphi_{1/\lambda}}f(\cdot,\xi)\|_{L^p}$
for each $0<\lambda\le 1$.
Hence,
by the same argument as in the proof of Lemma \ref{3.8},
we have
\begin{align*}
&\|V_{\varphi}(f_{\lambda})\|_{L^{p,\infty}}
=\lambda^{-n(1/p+1)}\|V_{\varphi_{1/\lambda}}f\|_{L^{p,\infty}}
\ge \lambda^{-n(1/p+1)}
\|V_{\varphi_{1/\lambda}}f(\cdot,0)\|_{L^p} \\
&\ge C\lambda^{-n(1/p+1)}
\left( \sum_{\ell \in \Z^n}
|\lambda^n \hat{\varphi}(\lambda \ell)|^p
\right)^{1/p}
\ge C\lambda^{-n/p}\left(
\sum_{\scriptstyle |\ell_j|\le 1/\lambda,
\atop \scriptstyle j=1,\cdots,n}
|\hat{\varphi}(\lambda \ell)|^p
\right)^{1/p}
\ge C\lambda^{-2n/p}
\end{align*}
for all $0<\lambda \le 1$.
Combining Lemma \ref{3.6} (1),
we get
$\|f\|_{M^{p,\infty}} \sim \lambda^{-2n/p}$
in the case $0<\lambda \le 1$.
The proof is complete.
\end{proof}
We are now ready to prove
Theorem \ref{1.1} (1) with $(1/p,1/q) \in I_3^*$
and (2) with $(1/p,1/q) \in I_3$.

\medskip
\noindent
{\it Proof of Theorem \ref{1.1} (2) with $(1/p,1/q) \in I_3$}.
We note that $\mu_2(p,q)=-2/p+1/q$ if $(1/p,1/q) \in I_3$.
Let $1 \le p \le 2$ and $1 \le q \le \infty$
be such that $(1/p,1/q) \in I_3$
and $1/p+1/q \ge 1$.
We note that,
if $(1/p,1/q) \in I_3$ and $1/p+1/q \ge 1$,
then $(1/p,1/q) \in I_2^*$ and $1/p \ge 1/q$.
Then, by Lemma \ref{3.2},
there exists a constant $C>0$ such that
\begin{equation}\label{(3.14)}
\|f_{\lambda}\|_{M^{p,q}}
\le C\lambda^{-n(2/p-1/q)}\|f\|_{M^{p,q}}
\quad \text{for all} \
f \in M^{p,q}(\R^n)
\ \text{and} \
0<\lambda \le 1.
\end{equation}
This is the first part of Theorem \ref{1.1} (2)
with $(1/p,1/q) \in I_3$
and $1/p+1/q \ge 1$.
Let $1 \le p \le 2$ and $2 \le q < \infty$
be such that $(1/p,1/q) \in I_3$
and $1/p+1/q \le 1$.
(\ref{(3.14)}) implies
\begin{equation}\label{(3.15)}
\|f_{\lambda}\|_{M^{r,r'}}
\le C\lambda^{-n(2/r-1/r')}\|f\|_{M^{r,r'}}
=C\lambda^{-n(3/r-1)}\|f\|_{M^{r,r'}}
\end{equation}
for all $f \in M^{r,r'}(\R^n)$
and $0<\lambda \le 1$,
where $1\le r \le 2$.
Take $1 \le r \le 2$ and $0\le \theta \le 1$
such that $1/p=(1-\theta)/2+\theta/r$
$1/q=(1-\theta)/\infty+\theta/r'$.
Since $q \neq \infty$,
we note that $r' \neq \infty$.
Then, by the interpolation theorem,
Lemma \ref{3.4} and (\ref{(3.15)})
give
\[
\|f_{\lambda}\|_{M^{p,q}}
\le C\left( \lambda^{-n} \right)^{1-\theta}
\left( \lambda^{-n(3/r-1)} \right)^{\theta}
\|f\|_{M^{p,q}}
\]
for all $f \in M^{p,q}(\R^n)$
and $0<\lambda \le 1$.
Since $1-\theta=2(1/p-\theta/r)$
and $\theta/r=\theta-1/q$,
we have
\begin{align*}
\|f_{\lambda}\|_{M^{p,q}}
&\le C\lambda^{-n(2(1/p-\theta/r)+3\theta/r-\theta)}
\|f\|_{M^{p,q}} \\
&=C\lambda^{-n(2/p+\theta/r-\theta)}
\|f\|_{M^{p,q}}
=C\lambda^{-n(2/p-1/q)}
\|f\|_{M^{p,q}}
\end{align*}
for all $f \in M^{p,q}(\R^n)$
and $0<\lambda \le 1$.
In the case $q=\infty$,
by Lemma \ref{3.6} (1),
we have nothing to prove.
Hence, we obtain the first part of Theorem \ref{1.1} (2)
with $(1/p,1/q) \in I_3$
and $1/p+1/q \le 1$.

By using Lemma \ref{3.8} (or \ref{3.9}),
we can prove the second part of Theorem \ref{1.1} (2)
with $(1/p,1/q) \in I_3$
in the same way
as in the proof of the second part of Theorem \ref{1.1} (2)
with $(1/p,1/q) \in I_2$.
We omit the proof.

\medskip
\noindent
{\it Proof of Theorem \ref{1.1} (1) with $(1/p,1/q) \in I_3^*$}.
We note that $\mu_1(p,q)=-2/p+1/q$ if $(1/p,1/q) \in I_3^*$.
In every case except for $(p,q)\neq (\infty,1)$,
by duality, Theorem \ref{1.1} (2) with $(1/p',1/q') \in I_3$
and the same argument as in the proof of Theorem \ref{1.1} (1)
with $(1/p,1/q) \in I_1^*$,
we can prove Theorem \ref{1.1} (1)
with $(1/p,1/q) \in I_3^*$.

For the first part of Theorem \ref{1.1} (1)
with $(p,q)=(\infty,1)$,
by (\ref{(3.11)}),
we have nothing to prove.
By using the interpolation theorem,
we can prove the second part in the same way
as in the proof of Theorem \ref{1.1} (1)
with $q=\infty$.
\section{The inclusion between Besov spaces and modulation spaces}
\label{section4}
In this section,
we prove Theorem \ref{1.2}
which appeared in the introduction.
It is sufficient to prove the first statement only because
the first one implies the second one by the duality argument and the
elementary relation
\[
\nu_2(p,q)=-\nu_1(p',q').
\]
See also Section 2 for the dual spaces of the modulation spaces
(Lemma \ref{2.4}) and Besov spaces.
\par
For the preparation
to prove Theorem \ref{1.2} (1)
with $(1/p,1/q) \in I_1^*$,
we show three lemmas in the below.
We denote by $B$ the tensor product of B-spline of degree $2$,
that is
\[
B(t)=\prod_{j=1}^n \chi_{[-1/2,1/2]}*\chi_{[-1/2,1/2]}(t_j),
\]
where $t=(t_1,\cdots, t_n) \in \R^n$.
We note that 
${\rm supp}\, B \subset [-1,1]^n$
and $\calF^{-1} B \in M^{p,q}(\R^n)$
for all $1 \le p,q \le \infty$.
\begin{lem}\label{4.1}
Let $1 \le p,q \le \infty$,
$(p,q)\neq (1,\infty), (\infty,1)$ and $\epsilon>0$.
Suppose that $\psi \in \calS(\R^n)$
satisfies
$\psi=1$ on $\{\xi: |\xi|\le 1/2\}$
and ${\rm supp}\, \psi \subset \{\xi: |\xi|\le 1\}$.
Set
\[
f(t)=\sum_{\ell \neq 0}
|\ell|^{-n/p-\epsilon}\, \Psi(t-\ell)
\quad \text{in} \quad \calS'(\R^n),
\]
where $\Psi=\calF^{-1}\psi$.
Then $f \in M^{p,q}(\R^n)$
and there exists a constant $C>0$ such that
$\|f_{\lambda}\|_{M^{p,q}}
\ge C\lambda^{-n/p-\epsilon}$
for all $\lambda \ge 2\sqrt{n}$.
\end{lem}
\begin{proof}
In the same way as in the proof of Lemma \ref{3.7},
we can prove $f \in M^{p,q}(\R^n)$.
We consider the second part.
Let $\lambda \ge 2\sqrt{n}$.
Since $\psi(\cdot/\lambda)=1$
on $[-1,1]^n$,
we have
\begin{align*}
&\int_{\R^n}
\Psi(\lambda t-\ell)\, (\calF^{-1}B)(t)\, dt
=(2\pi)^{-n}\lambda^{-n}
\int_{\R^n}e^{-i(\ell/\lambda)\cdot t}\,
\psi(t/\lambda)\, B(t)\, dt \\
&=(2\pi)^{-n}\lambda^{-n}
\int_{\R^n}e^{-i(\ell/\lambda)\cdot t}\, B(t)\, dt
=(2\pi)^{-n}\lambda^{-n}
\prod_{j=1}^n\left(
\frac{\sin \ell_j/2\lambda}{\ell_j/2\lambda} \right)^2.
\end{align*}
We note that
$\prod_{j=1}^n \left\{ (\sin \xi_j)/\xi_j \right\}^2 \ge C$
on $[-\pi/2,\pi/2]^n$
for some constant $C>0$.
Since $\calF^{-1}B \in M^{p',q'}(\R^n)$,
by Lemmas \ref{2.2} and \ref{2.4},
we get
\begin{align*}
&\|f_{\lambda}\|_{M^{p,q}}
=\sup_{\|g\|_{M^{p',q'}}=1}
\left| \langle f_{\lambda},g \rangle_M \right|
\ge \|\calF^{-1}B\|_{M^{p',q'}}^{-1}
\left| \langle f_{\lambda}, \calF^{-1}B \rangle_M \right| \\
&=\|\calF^{-1}B\|_{M^{p',q'}}^{-1}
\left| \sum_{\ell \neq 0}|\ell|^{-n/p-\epsilon}\,
\frac{1}{\|\Phi\|_{L^2}^2}
\int_{\R^{2n}}
V_{\Phi}[\Psi(\lambda \cdot-\ell)](x,\xi)\,
\overline{V_{\Phi}[\calF^{-1}B](x,\xi)}
\, dx \, d\xi \right| \\
&=\|\calF^{-1}B\|_{M^{p',q'}}^{-1}
\left| \sum_{\ell \neq 0}|\ell|^{-n/p-\epsilon}\,
\int_{\R^n}
\Psi(\lambda t-\ell)\, (\calF^{-1}B)(t)\, dt \right| \\
&=(2\pi)^{-n}\|\calF^{-1}B\|_{M^{p',q'}}^{-1}\lambda^{-n}
\left| \sum_{\ell \neq 0}|\ell|^{-n/p-\epsilon}\,
\prod_{j=1}^n \left( \frac{\sin \ell_j/2\lambda}{\ell_j/2\lambda}
\right)^2 \right| \\
&\ge (2\pi)^{-n}\|\calF^{-1}B\|_{M^{p',q'}}^{-1}\lambda^{-n}
\sum_{\scriptstyle 0 <|\ell_j| \le \lambda \pi,
\atop \scriptstyle j=1,\cdots,n}
|\ell|^{-n/p-\epsilon}\,
\prod_{j=1}^n
\left( \frac{\sin \ell_j/2\lambda}{\ell_j/2\lambda} \right)^2 \\
&\ge C\lambda^{-n}\lambda^{-n/p-\epsilon}
\sum_{\scriptstyle 0 <|\ell_j| \le \lambda \pi,
\atop \scriptstyle j=1,\cdots,n} 1
\ge C\lambda^{-n/p-\epsilon}
\end{align*}
for all $\lambda \ge 2\sqrt{n}$.
The proof is complete.
\end{proof}
\begin{lem}\label{4.2}
Suppose that $1 \le p,q \le \infty$,
$(p,q)\neq (1,\infty), (\infty,1)$ and $\epsilon>0$.
Let $\psi \in \calS(\R^n)$ be as in Lemma \ref{4.1}.
Set
\[
f(t)=e^{8it_1}\, \sum_{\ell \neq 0}
|\ell|^{-n/p-\epsilon}\, \Psi(t-\ell)
\quad \text{in} \quad \calS'(\R^n),
\]
where $t=(t_1,\cdots,t_n) \in \R^n$ and $\Psi=\calF^{-1}\psi$.
Then $f \in M^{p,q}(\R^n)$
and there exists a constant $C>0$ such that
$\|f_{\lambda}\|_{M^{p,q}}
\ge C\lambda^{-n/p-\epsilon}$
for all $\lambda \ge 2\sqrt{n}$.
\end{lem}
\begin{proof}
Let
$g(t)=\sum_{\ell \neq 0}|\ell|^{-n/p-\epsilon}\, \Psi(t-\ell)$.
Since $f=M_{8 e_1}g$ and $f_{\lambda}=M_{8\lambda e_1}g_{\lambda}$,
we have
$V_{\Phi}(f_{\lambda})(x,\xi)=V_{\Phi}(g_{\lambda})(x,\xi-8\lambda e_1)$,
where $e_1=(1,0,\cdots,0)$.
This gives
$\|f_{\lambda}\|_{M^{p,q}}=\|g_{\lambda}\|_{M^{p,q}}$.
Hence,
by Lemma \ref{4.1},
we obtain Lemma \ref{4.2}.
\end{proof}
\begin{lem}\label{4.3}
Suppose that $1 \le p,q \le \infty$,
$s \in \R$
and $\epsilon>0$.
Let $\psi \in \calS(\R^n)$ be as in Lemma \ref{4.1}.
Set
\[
f(t)=e^{8it_1}\,
\sum_{\ell \neq 0}
|\ell|^{-n/p-\epsilon}\, \Psi(t-\ell)
\quad \text{in} \quad \calS'(\R^n),
\]
where
$t=(t_1,\cdots,t_n) \in \R^n$ and $\Psi=\calF^{-1}\psi$.
Then there exists a constant $C>0$ such that
$\|f_{2^k}\|_{B_s^{p,q}}
\le C2^{k(s-n/p)}$
for all $k \in \Z_+$.
\end{lem}
\begin{proof}
Let $k \in \Z_+$.
Since ${\rm supp}\, \varphi_0 \subset \{\xi: |\xi| \le 2\}$,
${\rm supp}\, \varphi_j
\subset \{\xi: 2^{j-1} \le |\xi| \le 2^{j+1}\}$
if $j \ge 1$,
and
${\rm supp}\, \psi(\cdot/2^k-8e_1)$
$\subset$
$\{\xi:|\xi-2^{k+3}e_1|\le 2^k \}$,
we see that
\begin{align*}
&\int_{\R^n}
\Phi_j(x-t)\,
\left( e^{8 i (2^k t_1)}\, \Psi(2^k t -\ell) \right)\,
dt \\
&=(2\pi)^{-n}
\int_{\R^n}
e^{ix\cdot t}\,
\varphi_j(t)\,
\left( 2^{-kn}e^{-i\ell\cdot(t/2^k-8e_1)}\,
\psi(t/2^k-8e_1) \right) dt \\
&=
\begin{cases}
(2\pi)^{-n}e^{8i\ell_1}\,
\int_{\R^n}
e^{i(2^kx-\ell)\cdot t}\,
\varphi_j(2^k t)\, \psi(t-8e_1)\, dt,
&\text{if} \quad k+2 \le j \le k+4 \\
0,
&\text{otherwise}.
\end{cases}
\end{align*}
Hence,
\begin{align*}
&\left| \Phi_j*(f_{2^k})(x) \right|
=\left| \sum_{\ell \neq 0}
|\ell|^{-n/p-\epsilon}\,
\int_{\R^n}
\Phi_j(x-t)\,
\left( e^{8 i (2^k t_1)}\, \Psi(2^k t -\ell) \right)\,
dt \right| \\
&\le C\sum_{\ell \neq 0}
|\ell|^{-n/p-\epsilon}
\left| \int_{\R^n}
\left\{ (1+|2^kx-\ell|^2)^{-n}
(I-\Delta_t)^n e^{i(2^kx-\ell)\cdot t} \right\}
\varphi_j(2^k t)\, \psi(t-8e_1)\, dt \right| \\
&\le C\sum_{\ell \neq 0}
|\ell|^{-n/p-\epsilon}\, (1+|2^kx-\ell|^2)^{-n},
\end{align*}
where $k+2 \le j \le k+4$.
On the other hand,
$\Phi_j*(f_{2^k})=0$ if $j<k+2$ or $j>k+4$.
Thus,
$\|\Phi_j*(f_{2^k})\|_{L^p} \le C2^{-kn/p}$
if $k+2 \le j \le k+4$,
and $\|\Phi_j*(f_{2^k})\|_{L^p}=0$ if $j<k+2$ or $j>k+4$.
Therefore,
\[
\|f_{2^k}\|_{B_s^{p,q}}
=\left( \sum_{j=k+2}^{k+4}
2^{jsq}\|\Phi_j*(f_{2^k})\|_{L^p}^q \right)^{1/q}
\le C2^{-kn/p}\left( \sum_{j=k+2}^{k+4}2^{jsq} \right)^{1/q}
\le C2^{k(s-n/p)}.
\]
The proof is complete.
\end{proof}
We are now ready to prove Theorem \ref{1.2} (1)
with $(1/p,1/q) \in I_1^*$.

\medskip
\noindent
{\it Proof of Theorem \ref{1.2} (1) with $(1/p,1/q) \in I_1^*$}.
Let $(1/p,1/q) \in I_1^*$
and $(p,q)\neq (1,\infty)$.
Then $\nu_1(p,q)=0$.
We assume that
$B_s^{p,q}(\R^n) \hookrightarrow M^{p,q}(\R^n)$,
where $s<0$.
Set $s=-\epsilon$, where $\epsilon>0$.
For this $\epsilon$,
we define $f$ by
\[
f(t)=e^{8it_1}\,
\sum_{\ell \neq 0}
|\ell|^{-n/p-\epsilon/2}\, \Psi(t-\ell),
\]
where
$t=(t_1,\cdots,t_n) \in \R^n$, $\Psi=\calF^{-1}\psi$
and $\psi$ is as in Lemma \ref{4.1}.
Then,
by Lemmas \ref{4.2} and \ref{4.3},
we have
\[
C_1 2^{-k(n/p+\epsilon/2)}
\le \|f_{2^k}\|_{M^{p,q}}
\le C_2 \|f_{2^k}\|_{B_s^{p,q}}
\le C_3 2^{k(s-n/p)}
=C_3 2^{-k(n/p+\epsilon)}
\]
for any large integer $k$.
However, this is contradiction.
Hence, $s$ must satisfy $s \ge 0$.

We next consider the case $(p,q)=(1,\infty)$.
Assume that
$B_s^{1,\infty}(\R^n) \hookrightarrow M^{1,\infty}(\R^n)$.
Let $\psi \in \calS(\R^n)\setminus \{0\}$
be such that
${\rm supp}\, \psi \subset \{\xi: 1/2 \le |\xi| \le 2\}$.
Since $M^{1,\infty}(\R^n) \hookrightarrow \calF L^{\infty}(\R^n)$
(\cite[Proposition 1.7]{Toft}),
we see that
\[
2^{-kn}\|\psi\|_{L^{\infty}}
=\|\calF [\Psi_{2^k}]\|_{L^{\infty}}
\le C\|\Psi_{2^k}\|_{M^{1,\infty}}
\quad \text{for all} \ k \in \Z_+,
\]
where $\Psi=\calF^{-1}\psi$.
On the other hand,
it is easy to show that
\[
\|\Psi_{2^k}\|_{B_s^{1,\infty}}
\le C2^{k(s-n)}
\quad \text{for all} \ k \in \Z_+.
\]
Hence, by our assumption,
we get
\[
2^{-kn}\|\psi\|_{L^{\infty}}
\le C_1\|\Psi_{2^k}\|_{M^{1,\infty}}
\le C_2\|\Psi_{2^k}\|_{B_s^{1,\infty}}
\le C_3 2^{k(s-n)}
\]
for all $k \in \Z_+$.
This implies $s \ge 0$.
The proof is complete.

\medskip
Our next goal is to prove Theorem \ref{1.2} (1)
with $(1/p,1/q)\in I_2^*$.
We remark the following fact,
and give the proof for reader's convenience.
\begin{lem}[{\cite[Proposition 1.1]{Sugimoto}}]\label{4.4}
Let $1\le p,q \le \infty$ and $s>0$.
Then there exists a constant $C>0$ such that
\[
\|f_{\lambda}\|_{B_s^{p,q}}
\le C\lambda^{s-n/p}\|f\|_{B_s^{p,q}}
\quad \text{for all} \
f \in B_s^{p,q}(\R^n)
\ \text{and} \ \lambda \ge 1.
\]
\end{lem}
\begin{proof}
Let $j_0 \in \Z_+$ be such that
$2^{j_0} \le \lambda <2^{j_0+1}$.
Since $\sum_{j=0}^{\infty}\varphi_j(\xi)=1$
for all $\xi \in \R^n$,
we see that
\[
\varphi_j(\lambda \xi)
=\sum_{\ell=-2}^{1}
\varphi_j(\lambda \xi)\,
\varphi_{j+\ell}(2^{j_0}\xi)
\quad \text{for all} \
\xi \in \R^n
\ \text{and} \
j \in \Z_+,
\]
where $\varphi_{j+\ell}=0$ if $j+\ell<0$.
Hence, by Young's inequality,
we have
\begin{align*}
\|f_{\lambda}\|_{B_s^{p,q}}
&=\left( \sum_{j=0}^{\infty}2^{jsq}
\|\calF^{-1}[\varphi_j \, \widehat{f_{\lambda}}]\|_{L^p}^q
\right)^{1/q}
=\lambda^{-n/p}
\left( \sum_{j=0}^{\infty}2^{jsq}
\|\calF^{-1}[\varphi_j(\lambda \cdot) \, \hat{f}]\|_{L^p}^q
\right)^{1/q} \\
&\le \lambda^{-n/p}\sum_{\ell=-2}^{1}
\left( \sum_{j=0}^{\infty}2^{jsq}
\|\calF^{-1}[\varphi_j(\lambda \cdot) \,
\varphi_{j+\ell}(2^{j_0}\cdot) \, \hat{f}]\|_{L^p}^q
\right)^{1/q} \\
&\le \lambda^{-n/p}\sum_{\ell=-2}^{1}
\left\{ \sum_{j=0}^{\infty}2^{jsq}
\left( \|\calF^{-1}[\varphi_j(\lambda \cdot)]\|_{L^1}
\|\calF^{-1}[\varphi_{j+\ell}(2^{j_0}\cdot) \, \hat{f}]\|_{L^p}
\right)^q \right\}^{1/q} \\
&\le C\lambda^{-n/p}
\left( \sum_{j=0}^{\infty}2^{jsq}
\|\calF^{-1}[\varphi_{j}(2^{j_0}\cdot) \, \hat{f}]\|_{L^p}^q
\right)^{1/q} \\
&=C\lambda^{-n/p}
\left\{ \left(\sum_{j=0}^{j_0}+\sum_{j=j_0+1}^{\infty}
\right) 2^{jsq}
\|\calF^{-1}[\varphi_{j}(2^{j_0}\cdot) \, \hat{f}]\|_{L^p}^q
\right\}^{1/q}.
\end{align*}
For the first term,
we see that
\begin{align*}
&\sum_{j=0}^{j_0}2^{jsq}
\|\calF^{-1}[\varphi_{j}(2^{j_0}\cdot) \, \hat{f}]\|_{L^p}^q
=\sum_{j=0}^{j_0}2^{jsq}
\|\calF^{-1}[\varphi_{j}(2^{j_0}\cdot) \,
(\varphi_0+\varphi_1+\varphi_2)\, \hat{f}]\|_{L^p}^q \\
&\le C\sum_{j=0}^{j_0}2^{jsq}
\|\calF^{-1}[(\varphi_0+\varphi_1+\varphi_2)\, \hat{f}]\|_{L^p}^q
\le C\left( 2^{j_0 s}\|f\|_{B_s^{p,q}}\right)^q
\le C\left( \lambda^s\|f\|_{B_s^{p,q}}\right)^q.
\end{align*}
For the second term,
we have
\[
\sum_{j=j_0+1}^{\infty}2^{jsq}
\|\calF^{-1}[\varphi_{j}(2^{j_0}\cdot) \, \hat{f}]\|_{L^p}^q
=\sum_{j=j_0+1}^{\infty}2^{jsq}
\|\calF^{-1}[\varphi_{j-j_0}\, \hat{f}]\|_{L^p}^q
\le \left( \lambda^{s}\|f\|_{B_s^{p,q}}\right)^q.
\]
Combining these estimates,
we obtain the desired result.
\end{proof}
We are now ready to prove Theorem \ref{1.2} (1)
with $(1/p,1/q) \in I_2^*$.

\medskip
\noindent
{\it Proof of Theorem \ref{1.2} (1) with $(1/p,1/q) \in I_2^*$}.
Let $(1/p,1/q) \in I_2^*$.
Then $\nu_1(p,q)=1/p+1/q-1$.
If $(1/p,1/q) \in I_2^*$ and $1/p+1/q=1$
then $(1/p,1/q) \in I_1^*$,
and we have already proved this case
in Theorem \ref{1.2} (1) with $(1/p,1/q) \in I_1^*$.
Hence, we may assume $1/p+1/q>1$.
Suppose that
$B_s^{p,q}(\R^n) \hookrightarrow M^{p,q}(\R^n)$,
where $s<n(1/p+1/q-1)$.
Then, since $n(1/p+1/q-1)>0$,
we can take $s_0>0$ such that
$s\le s_0 <n(1/p+1/q-1)$.
Let $\varphi$ be the Gauss function.
By Lemma \ref{2.1},
we see that
$\|\varphi_{\lambda}\|_{M^{p,q}}\ge C \lambda^{n(1/q-1)}$
for all $\lambda \ge 1$.
On the other hand,
by Lemma \ref{4.4},
we have
\[
\|\varphi_{\lambda}\|_{B_{s_0}^{p,q}}
\le C\lambda^{s_0-n/p}\|\varphi\|_{B_{s_0}^{p,q}}
\quad \text{for all} \ \lambda \ge 1.
\]
Hence,
using
$B_{s_0}^{p,q}(\R^n)
\hookrightarrow B_s^{p,q}(\R^n)
\hookrightarrow M^{p,q}(\R^n)$,
we get
\[
C_1 \lambda^{n(1/q-1)}
\le \|\varphi_{\lambda}\|_{M^{p,q}}
\le C_2 \|\varphi_{\lambda}\|_{B_{s_0}^{p,q}}
\le C_3 \lambda^{s_0-n/p}\|\varphi\|_{B_{s_0}^{p,q}}
\]
for all $\lambda \ge 1$.
However,
since $s_0-n/p<n(1/q-1)$,
this is contradiction.
Therefore, $s$ must satisfy $s \ge n(1/p+1/q-1)$.
The proof is complete.

\medskip
Our next goal is to prove Theorem \ref{1.2} (1)
with $(1/p,1/q) \in I_3^*$.
\begin{lem}\label{4.5}
Let $1 \le p\le \infty$,
$1\le q <\infty$ and $\epsilon >0$.
Suppose that
$\varphi, \psi \in \calS(\R^n)\setminus \{0\}$ satisfy
${\rm supp}\, \varphi \subset [-1/8,1/8]^n$,
${\rm supp}\, \psi \subset [-1/2,1/2]^n$
and $\psi=1$ on $[-1/4,1/4]^n$.
For $j \in \Z_+$,
set
\[
f^j(t)=2^{-jn/p}
\sum_{\scriptstyle 0<|k_j|\le 2^j, \atop \scriptstyle j=1,\cdots , n}
|k|^{-n/p-\epsilon}\,
e^{ik\cdot t/2^j}\, \Psi(t/2^j-k),
\]
where $\Psi=\calF^{-1}\psi$.
Then $f^j \in M^{p,q}(\R^n)$
and there exists a constant $C>0$ such that
\[
\| V_{\Phi}[(f^j)_{2^j}] \|_{L^{p,q}}
\ge C2^{-jn(2/p-1/q)-j\epsilon}
\quad \text{for all} \ j \in \Z_+,
\]
where
$\Phi=\calF^{-1}\varphi$.
\end{lem}
\begin{proof}
Since $f^j \in \calS(\R^n)$,
we have $f^j \in M^{p,q}(\R^n)$.
We consider the second part.
Note that
${\rm supp}\, \varphi(\cdot-\xi) \subset \ell + [-1/4,1/4]^n$
for all $\ell \in \Z^n$
and $\xi \in \ell+[-1/8,1,8]^n$.
Since
${\rm supp}\, \psi(\cdot-k) \subset k+[-1/2,1/2]^n$
and $\psi(t-k)=1$ if $t \in k+[-1/4, 1/4]^n$,
it follows that
\begin{align*}
&\|V_{\Phi}[(f^j)_{2^j}]\|_{L^{p,q}}
\ge \bigg\{
\sum_{\ell \in \Z^n}
\int_{\ell +[-1/8,1/8]^n} \\ 
&\times \bigg( \int_{\R^n} \bigg| 2^{-jn/p}
\sum_{\scriptstyle 0<|k_j|\le 2^j, \atop \scriptstyle j=1,\cdots , n}
|k|^{-n/p-\epsilon}
\int_{\R^n}e^{ik\cdot t}\, \Psi(t-k)\,
\overline{\Phi(t-x)}\, e^{-i\xi \cdot t}\, dt
\bigg|^p dx \bigg)^{q/p} d\xi \bigg\}^{1/q} \\
&\ge(2\pi)^{-n}2^{-jn/p}\bigg\{
\sum_{\scriptstyle 0<|\ell_j|\le 2^j, \atop \scriptstyle j=1,\cdots , n}
\int_{\ell +[-1/8,1/8]^n} \\ 
&\quad \times \bigg( \int_{\R^n} \bigg|
\sum_{\scriptstyle 0<|k_j|\le 2^j, \atop \scriptstyle j=1,\cdots , n}
|k|^{-n/p-\epsilon}e^{i|k|^2}
\int_{\R^n}e^{-ik\cdot t}\, \psi(t-k)\,
\overline{\varphi(t-\xi)}\, e^{ix\cdot t}\, dt
\bigg|^p dx \bigg)^{q/p} d\xi \bigg\}^{1/q} \\
&=(2\pi)^{-n}2^{-jn/p}\bigg\{
\sum_{\scriptstyle 0<|\ell_j|\le 2^j, \atop \scriptstyle j=1,\cdots , n}
\int_{\ell +[-1/8,1/8]^n} \\
&\quad \times
\bigg( \int_{\R^n} \bigg|
|\ell|^{-n/p-\epsilon}
\int_{\R^n}e^{i(x-\ell)\cdot t}\,
\overline{\varphi(t-\xi)}\, dt
\bigg|^p dx \bigg)^{q/p} d\xi \bigg\}^{1/q} \\
&=2^{-jn/p}\bigg\{
\sum_{\scriptstyle 0<|\ell_j|\le 2^j, \atop \scriptstyle j=1,\cdots , n}
|\ell|^{-(n/p+\epsilon)q}
\int_{\ell +[-1/8,1/8]^n} 
\| \Phi(-\cdot+\ell)\|_{L^p}^q\,
d\xi \bigg\}^{1/q} \\
&=4^{-n/q}\|\Phi\|_{L^p}2^{-jn/p}\bigg\{
\sum_{\scriptstyle 0<|\ell_j|\le 2^j, \atop \scriptstyle j=1,\cdots , n}
|\ell|^{-(n/p+\epsilon)q}\bigg\}^{1/q} \\
&\ge C_n2^{-jn/p}
2^{-j(n/p+\epsilon)}\bigg\{
\sum_{\scriptstyle 0<|\ell_j|\le 2^j, \atop \scriptstyle j=1,\cdots , n}
1 \bigg\}^{1/q}
\ge C_n 2^{-jn(2/p-1/q)-j\epsilon}
\end{align*}
for all $j \in \Z_+$. The proof is complete.
\end{proof}
\begin{lem}\label{4.6}
Suppose that $1\le p,q \le \infty$ and $s>0$.
Let $f^j$ be as in Lemma \ref{4.5}.
Then there exists a constant $C>0$ such that
$\|(f^j)_{2^j}\|_{B_s^{p,q}} \le C2^{j(s-n/p)}$
for all $j \in \Z_+$.
\end{lem}
\begin{proof}
By Lemma \ref{4.4},
we have
$\|(f^j)_{2^j}\|_{B_s^{p,q}} \le C2^{j(s-n/p)}\|f^j\|_{B_s^{p,q}}$
for all $j \in \Z_+$.
Hence,
it is enough to prove that
$\sup_{j \in \Z_+}\|f^j\|_{B_s^{p,q}}<\infty$.
Since
\[
\widehat{f^j}(\xi)=2^{jn(1-1/p)}
\sum_{\scriptstyle 0<|k_j|\le 2^j, \atop \scriptstyle j=1,\cdots , n}
|k|^{-n/p-\epsilon}\,
e^{-ik\cdot (2^j\xi-k)}\, \psi(2^j\xi-k)
\]
and
${\rm supp}\, \psi(2^j\cdot-k) \subset k/2^j+[-2^{-(j+1)},2^{-(j+1)}]^n$,
we see that
${\rm supp}\, \widehat{f^j} \subset \{\xi:|\xi| \le 2\sqrt{n}\}$.
Let $\ell_0$ be such that $2^{\ell_0-1} \ge 2\sqrt{n}$.
Then,
\begin{align*}
\|f^j\|_{B_s^{p,q}}
&=\left(\sum_{\ell=0}^{\ell_0-1}
2^{\ell sq}\|\Phi_{\ell}*f^j\|_{L^p}^q \right)^{1/q}
\le \left(\sum_{\ell=0}^{\ell_0-1}
2^{\ell sq}(\|\Phi_{\ell}\|_{L^1}\|f^j\|_{L^p})^q \right)^{1/q}
=C_n \|f^j\|_{L^p}.
\end{align*}
Therefore, it is enough to show that
$\sup_{j \in \Z_+}\|f^j\|_{L^p}<\infty$.
By a change of variable,
we have
\begin{align*}
\|f^j\|_{L^p}
&=\bigg( \int_{\R^n} \bigg|
\sum_{\scriptstyle 0<|k_j|\le 2^j, \atop \scriptstyle j=1,\cdots , n}
|k|^{-n/p-\epsilon}\, e^{ik\cdot t}\,
\Psi(t-k) \bigg|^p \, dt \bigg)^{1/p} \\
&\le \bigg\{ \sum_{m \in \Z^n}
\int_{m+[-1/2,1/2]^n} \bigg(
\sum_{k \neq 0}|k|^{-n/p-\epsilon}\, |\Psi(t-k)|
\bigg)^p \, dt \bigg\}^{1/p} \\
&\le C\bigg\{ \sum_{m \in \Z^n} \bigg(
\sum_{k \neq 0}|k|^{-n/p-\epsilon}\, (1+|m-k|)^{-n-1}
\bigg)^p \bigg\}^{1/p}<\infty
\end{align*}
for all $j \in \Z_+$.
The proof is complete.
\end{proof}
We are now ready to prove Theorem \ref{1.2} (1)
with $(1/p,1/q) \in I_3^*$.

\medskip
\noindent
{\it Proof of Theorem \ref{1.2} (1) with $(1/p,1/q) \in I_3^*$}.
Let $(1/p,1/q) \in I_3^*$.
Then $\nu_1(p,q)=-1/p+1/q$.
If $(1/p,1/q) \in I_3^*$ and $p=q$
then $(1/p,1/q) \in I_1^*$,
and we have already proved this case
in Theorem \ref{1.2} (1) with $(1/p,1/q) \in I_1^*$.
Hence, we may assume $1/q>1/p$.
Note that $q \neq \infty$.
Suppose that $B_s^{p,q}(\R^n) \hookrightarrow M^{p,q}(\R^n)$,
where $s < -n(1/p-1/q)$.
Then, since $-n(1/p-1/q)>0$,
we can take $s_0>0$ such that $s \le s_0 < -n(1/p-1/q)$.
Set $s_0=-n(1/p-1/q)-\epsilon$,
where $\epsilon>0$.
For this $\epsilon$,
we define $f^j$ by
\[
f^j(t)=2^{-jn/p}
\sum_{\scriptstyle 0<|k_j|\le 2^j, \atop \scriptstyle j=1,\cdots , n}
|k|^{-n/p-\epsilon/2}\,
e^{ik\cdot t/2^j}\, \Psi(t/2^j-k),
\]
where $j \in \Z_+$,
$\Psi=\calF^{-1}\psi$ and $\psi$ is as in Lemma \ref{4.5}.
Then,
since
$B_{s_0}^{p,q}(\R^n)
\hookrightarrow B_s^{p,q}(\R^n)
\hookrightarrow M^{p,q}(\R^n)$,
by Lemmas \ref{4.5} and \ref{4.6},
we get
\begin{align*}
C_1 2^{-jn(2/p-1/q)-j\epsilon/2}
&\le \|V_{\Phi}[(f^j)_{2^j}]\|_{L^{p,q}}
\le C_2\|(f^j)_{2^j}\|_{M^{p,q}} \\
&\le C_3\|(f^j)_{2^j}\|_{B_{s_0}^{p,q}}
\le C_4 2^{j(s_0-n/p)}
=C_4 2^{-jn(2/p-1/q)-j\epsilon}
\end{align*}
for all $j \in \Z_+$,
where $\Phi=\calF^{-1}\varphi$ and $\varphi$ is as in Lemma \ref{4.5}.
However,
this is contradiction.
Therefore, $s$ must satisfy $s \ge -n(1/p-1/q)$.
The proof is complete.

\end{document}